\documentclass[runningheads]{llncs}
\usepackage{graphicx}
\usepackage{amssymb}
\usepackage{amstext}
\usepackage{amsmath}
\usepackage{color,tabularx,multirow,hyperref}

\newcommand{\RR}{{\mathbb R}}

\begin{document}
\title{Optimal reachability and grasping\\ for a soft manipulator}
\author{Simone Cacace\inst{1}\and
Anna Chiara Lai\inst{2}\orcidID{0000-0003-2096-6753} \and
Paola Loreti\inst{2}}
\authorrunning{S. Cacace et al.}
\institute{Dipartimento di Matematica e Fisica, Universit\`{a} degli studi Roma Tre, Rome, Italy \\
\email{cacace@mat.uniroma3.it} \and
Dipartimento di Scienze di Base e Applicate per l'Ingegneria, Sapienza Universit\`{a} di Roma, Rome, Italy\\
\email{\{anna.lai,paola.loreti\}@sbai.uniroma1.it}}
\maketitle              % typeset the header of the contribution
\begin{abstract}
We investigate optimal reachability and grasping problems for a planar soft manipulator, from both a theoretical and numerical point of view. The underlying control model describes the evolution of the symmetry axis of the device, which is subject to inextensibility and curvature constraints, a bending moment and a curvature control. Optimal control strategies are characterized with tools coming from the optimal control theory of PDEs. We run some numerical tests in order to validate the model and to synthetize optimal control strategies. 

\end{abstract}
\section{\uppercase{Introduction}}
\label{sec:introduction}
In this paper we investigate reachability and grasping problems for planar soft manipulators in the framework of optimal control theory. The manipulator we are modeling is a three-dimensional body with an axial symmetry, non-uniform thickness and a fixed endpoint. The device is characterized by a structural resistance to bending, the \emph{bending moment}, and its soft material allows for elastic distorsions only below a fixed threshold, modelled via a \emph{bending constraint}. Moreover, the bending can be enforced by a pointwise \emph{bending control}, an internal angular force representing the control term of the system. Finally, we assume an \emph{inextensiblity constraint} on the manipulator: its structure does not allow for longitudinal stretching. Using the key morphological assumption of axial symmetry, we restrict our investigation to the evolution of the symmetry axis of the manipulator, modeled as a planar curve $q(s,t):[0,1]\times[0,T]\to\RR^2$, where $s$ is an arclength coordinate and $T>0$ is the final time. Physically, $q$ is modelled as an inextensible beam, whose mass represents the mass of the whole manipulator. The evolution is determined by  suitable reaction and friction forces, encompassed in a nonlinear, fourth order system of evolutive, controlled PDEs, generalizing the Euler-Bernoulli beam equation. In particular, the aforementioned bending constraints and controls of the manipulator on its symmetry axis are encoded as curvature constraints and controls, enforced by angular reaction forces. 

The optimal control problems that we address are formulated as a constrained minimization of appropriate cost functionals, involving a quadratic cost on the controls, see \eqref{js}, \eqref{dynamicfunc} for the reachability problem and \eqref{contactfunctional} for the grasping problem. 

The control model underlying our study and a first investigation of the related optimal reachability strategies were originally developed in \cite{CLL18}, then extended  in \cite{CLL19} to the case in which only a portion of the device is controlled. Here, beside a deeper investigation of the model and of its connection with hyper-redundant manipulators, we address the important problem of grasping a prescribed object. In our investigation the grasping problem is studied (and numerically solved) in a stationary setting. We point out that our approach, involving the theory of optimal control of PDEs, allows in a quite natural way for an extension to fully dynamic optimization: here we present the key, spatial component of the numerical solution of this challenging problem, while postponing its evolutive counterpart to future works.  

% Soft robotics is a fascinating research field conjugating countless applications --for instance  manipulation in constrained envinroments and human-robot interactions  \cite{nature,laschi1}-- and as many  theoretical and practical interesting open problems. 
It is almost impossible to give an exhaustive state of the art even on the particular framework of our paper, i.e., the motion planning of highly articulated and soft manipulators. Then we limit ourselves to refer to  the papers that mostly inspired our work: the seminal papers \cite{hyper,hyper2} were hyper-redundant manipulator were originally introduced, the paper \cite{kinematicssoft} for the  kinematic aspects and \cite{octopusrobotprescribed,dynamicoctopusrobot} for a discrete dynamical model. The parameter setting of our tests is motivated by the number theoretical approach presented in \cite{hand,icinco1,icinco2}. Finally, we refer to the book \cite{optimalpdebook} for an introduction to our theoretical background, that is, optimal control theory of PDEs.

The paper is organized as follows. Section \ref{model} is devoted to an overview of the model under exam, and Section \ref{control} to static and dynamic optimal reachability problems. In Section \ref{grasping} we investigate an optimal grasping problem for the stationary case. Finally in Section \ref{conclusions} we draw our conclusions.

\section{\uppercase{A control model for planar manipulators}}\label{model}
We revise, from a robotic perspective, the main steps of the modeling procedure in \cite{CLL18}. In particular, we provide an analytical description of the constraints as formal limit of a discrete particle system, which in turn models an hyper-redundant manipulator.
\subsection{From a discrete hyper-redundant to a soft manipulator}
Hyper-redundant manipulator are characterized by a very large, possibly infinite, number of actuable degrees of freedom: their investigation in the framework of soft robotics is motivated by the fact that they can be viewed as discretizations of continuum robots. Following this idea, we consider here a device composed by $N$ links and $N+1$ joints, whose position in the plane is given by the array $(q_0,\dots,q_N)$. We denote by $m_k$ the mass of the $k$-th joint, with $k=0,\dots,N$, and we consider negligible the mass of the corresponding links.
Moreover, given vectors $v_1, v_2\in\RR^2$, we assume standard notations for the Euclidean norm and the dot product, respectively $|v_1|$ and $v_1\cdot v_2$, and we define $v_1\times v_2:=v_1\cdot v_2^\bot$, where $v_2^\bot$ denotes the clockwise orthogonal vector to $v_2$. Finally, the positive part function is denoted by $(\cdot)_+$.\\

The manipulator is endowed with the following physical features:
\begin{itemize}
    \item \textbf{Inextensibility constraint.} The links are rigid and their lengths are given by $\ell$, so that the discrete counterpart of the inextensibility constraint is 
\begin{equation}
    |q_{k}-q_{k-1}|=\ell \qquad k=1,\dots,N.
\end{equation}
This constraint is imposed exactly, by defining for $k=1,...,N$
$$F_k(q,\sigma):=\sigma_k\left(|q_{k}-q_{k-1}|^2-\ell^2\right)\,,$$
where $\sigma_k$ is a Langrange multiplier.
 \item \textbf{Curvature constraint.} We assume that two consecutive links, say the $k$-th and the $k+1$-th, cannot form an angle larger than a fixed threshold $\alpha_k$. This can be formalized by the following condition on the joints:
$$(q_{k+1}-q_k)\cdot(q_{k}-q_{k-1})\ge \ell^2\,\cos(\alpha_k)\,.$$
This constraint is imposed via penalty method. In particular,
for $k=1,...,N$ we consider the elastic potential: 
$$G_k(q):=\nu_k\left(\cos(\alpha_k)-\frac{1}{\ell^2}(q_{k+1}-q_{k})\cdot(q_{k}-q_{k-1})\right)_+^2\,,$$
 where $\nu_k\geq0$ is a penalty parameter playing the role of an elastic constant,. 
\item\textbf{Bending moment.} Modeling an intrinsic resistance to leave the
position at rest, corresponding to null relative angles, is described by the following equality constraint:
$$(q_{k+1}-q_k)\times(q_{k}-q_{k-1})=0\,.$$
The related elastic potential with penalty parameter $\varepsilon_k>0$ is
$$B_k(q):=\varepsilon_k\Big( (q_{k+1}-q_{k})\times(q_{k}-q_{k-1})\Big)^2\,,$$

\item \textbf{Curvature control.} A curvature control in a discrete setting equals to impose the exact angle between the joints $q_{k-1}, q_k, q_{k+1}$, i.e., the  following equality constraint: 
$$(q_{k+1}-q_k)\times(q_{k}-q_{k-1})=\ell^2\,\sin(\alpha_k u_k)\,,$$
where $u_k\in[-1,1]$ is the control term. Note that the control set $[-1,1]$ is chosen in order to be consistent with the curvature constraint. Also in this case we enforce the constraint by penalty method, by considering:
$$H_k(q):=\mu_k\left(\sin(\alpha_ku_k)-\frac{1}{\ell^2}(q_{k+1}-q_{k})\times(q_{k}-q_{k-1})\right)^2\,,$$
where $\mu_k\geq0$ is a penalty parameter. Note that to set $\mu_k=0$ corresponds to deactivate the control of the $k$-th joint and let it evolve according the above constraints only. 
\end{itemize}
Note that the definition of $G_{k},B_k$ and $H_k$ in the cases $k=0$ and $k=N$ is made consistent by considering two ghost joints $q_{-1}:=q_0+\ell(0,1)$ and $q_{N+1}:=q_N+(q_{N}-q_{N-1})$ at the endpoints. 
\medskip

We now are in position to build the Lagrangian associated to the hyper-redundant manipulator, composed by a kinetic energy term and suitably rescaled elastic potentials:
$$\mathcal L_N(q,\dot q,\sigma):=\sum_{k=0}^N \frac{1}{2}m_k|\dot q_k|^2-\frac{1}{2\ell}F_k(q,\sigma)-\frac{1}{\ell^3}G_k(q)-\frac{1}{2\ell^5}B_k(q)-\frac{1}{2\ell}H_k(q).$$

% NESSUNA INFO SULLA DENSITA' $\rho_k$.

% we take the masses and the maximum bending angles 
% of the form $m_k=\ell\rho_k$ and $\alpha_k=\ell\omega_k$ respectively, for given mass distributions $\rho_k$ and curvatures $\omega_k$. 

Assume now that the hyper-redundant manipulator is indeed a discretization of our continuous manipulator, that is, there exist smooth functions $\nu, \mu, \varepsilon, \rho, \omega:[0,1]\to\RR^+$ and, for $T>0$, smooth functions 
$q:[0,1]\times[0,T]\to\RR^2$, $\sigma:[0,1]\times[0,T]\to\RR$ and $u:[0,1]\times[0,T]\to[-1,1]$ such that, for all $N$, $k=1,...,N$ and $t\in[0,T]$
$$
\nu_k=\nu(k\ell),\quad \mu_k=\mu(k\ell),\quad \varepsilon_k=\varepsilon(k\ell),$$
$$ \rho_k:=m_k/\ell=\rho(k\ell),\quad \omega_k:=\alpha_k/\ell=\omega(k\ell),
$$
$$
q_k(t)=q(k\ell,t),\quad \sigma_k(t)=\sigma(k\ell,t),\quad u_k(t)=u(k\ell,t)\,.
$$
 
Fixing the total length of the manipulator equal to $1$, so that $\ell=1/N$, we define the Lagrangian associated to the soft, continuous manipulator by the formal limit (see \cite{CLL18} for details)
$$\lim_{N\to+\infty}\mathcal L_N(q,\dot q, \sigma)=\mathcal L(q,\dot q, \sigma).$$
where
\begin{equation}\label{continuouslagrangian}
 \begin{split}\mathcal{L}(q,\sigma):=\int_0^1\Big(&\frac12\rho|q_t|^2-\frac12 \sigma(|q_s|^2 - 1)-\frac14\nu\left(|q_{ss}|^2-\omega^2\right)_+^2
 -\frac12\varepsilon|q_{ss}|^2
 \\-&\frac12\mu\left(\omega u-q_s\times q_{ss}\right)^2
 \Big)ds\end{split}
 \end{equation}
  and $q_t$, $q_{s}$, $q_{ss}$ denote partial derivatives in time and space respectively.
 
A comparison between discrete and continuous constraints and related potentials is displayed in Table \ref{constraints} and Table \ref{potentials}, respectively. 
  
\begin{table}
\centering\caption{\label{constraints} Exact constraint equations in both discrete and continuous settings.}
{\renewcommand\arraystretch{1.3}\begin{tabular}{|l|l|l|}\hline
 	{\bf Constraint}& \multicolumn{2}{c|}{\bf Constraint equation} \\\hline
 &Discrete&Continuous\\\cline{2-3}
  Inextensibility& $|q_k-q_{k-1}|=\ell$&$|q_s|=1$\\
Curvature &$(q_{k+1}-q_k)\cdot(q_{k}-q_{k-1})\ge \ell^2\,\cos(\alpha_k)$ &$|q_{ss}|\leq \omega$\\
Bending moment&$(q_{k+1}-q_k)\times(q_{k}-q_{k-1})=0\,$&$|q_{ss}|=0$\\
 Control &$(q_{k+1}-q_k)\times(q_{k}-q_{k-1})=\ell^2\,\sin(\alpha_k u_k)$&$q_s \times q_{ss}=\omega u$\\\hline
 \end{tabular}}\end{table}
 
 \begin{table}\caption{\label{potentials}  Potentials derived from penalty method in both discrete and continuous settings. The functions $\nu$ and $\mu$ represent non-uniform elastic constants. }
 {\renewcommand\arraystretch{1.3}\begin{tabular}{|l|l|l|}\hline
 	{\bf Constraint}& 
 	 \multicolumn{2}{c|} {\bf Penalization elastic potential}\\\hline
 &Discrete&Continuous\\\cline{2-3}
  Inextensibility&None&None\\
 Curvature &$\nu_k\left(\cos(\alpha_k)-\frac{1}{\ell^2}(q_{k+1}-q_{k})\cdot(q_{k}-q_{k-1})\right)_+^2$& $\nu(|q_{ss}|^2-\omega^2)^2_+$\\
Bending moment&$\varepsilon_k\Big( (q_{k+1}-q_{k})\times(q_{k}-q_{k-1})\Big)^2$&$\varepsilon |q_{ss}|^2$\\
 Control & $\mu_k\left(\sin(\alpha_ku_k)-\frac{1}{\ell^2}(q_{k+1}-q_{k})\times(q_{k}-q_{k-1})\right)^2$& $\mu\left(\omega u-q_s \times q_{ss}\right)^2$
 \\\hline
 \end{tabular}}
\end{table}
\subsection{Equations of motion}
Equations of motion for both the discrete and the continuous models can be derived applying the least action principle to the corresponding Lagrangians. 
Here we report only the continuous case, which is the building block for the optimal control problems we address in the next sections. Taking into account also some friction forces, we obtain the following 
system of nonlinear, evolutive, fourth order PDEs -- see \cite{CLL18} for details:
\begin{equation}\label{sysintro}
\left\{
\begin{array}{l}
\rho q_{tt}=\left(\sigma q_s-H q_{ss}^\bot\right)_s -\left(G q_{ss}+H q_{s}^\bot\right)_{ss}-\beta q_t-\gamma q_{sssst}\,.
\\
|q_s|^2=1 
\end{array}
\right.
\end{equation}
for $(s,t)\in(0,1)\times(0,T)$. We remark that the map $G:=G[q,\nu,\varepsilon,\omega]=\varepsilon+\nu\left(|q_{ss}|^2-\omega^2\right)_+$ encodes the bending moment and the curvature constraint, while the map $H:=H[q,\mu,u,\omega]=\mu\left(\omega u-q_s \times q_{ss}\right)$ corresponds to  the control term. The term  $-\beta q_t$ represents  an environmental viscous friction proportional to the velocity; the term $-\gamma q_{sssst}$ an internal viscous friction, proportional to the change in time of the curvature.
The system is completed with suitable initial data and with the following 
boundary conditions for $t\in(0,T)$:
 \begin{equation}\label{bond1}\begin{cases}
q(0,t)=(0,0),&\text{ (anchor point)}\\
q_s(0,t)=-(0,1)&\text{ (fixed tangent)}\\
q_{ss}(1,t)=0 &\text{ (zero bending moment)} \\
q_{sss}(1,t)=0 &\text{ (zero shear stress)}\\
\sigma(1,t)=0 &\text{ (zero tension)}.
\end{cases}\end{equation}
Note that the first two conditions are a modeling choice, whereas the free endpoint conditions emerge from the stationarity of the Lagrangian $\mathcal L$. 
 
 For reader's convenience, we summarize the model parameters in Table \ref{parameters}.
 \begin{table}[!h] \caption{A summary of the quantities and the functions involved in the soft manipulator control equation \eqref{sysintro}.}
    \label{parameters}
    \centering
   {\renewcommand\arraystretch{1.3} \begin{tabular}{|c|l|l|}
    \hline
    {\,\,\bf Name\,\,}&{\bf Description}&{\bf Type}
    \\\hline
    $q$&Axis parametrization&\multirow{2}{*}{Unknown of the equation}\\
    $\sigma$& Inextensibility constraint multiplier&\\
        \hline
$u$& Curvature control & Control\\
\hline
    $\rho$& Mass distribution& \multirow{2}{*}{Physical parameter}\\
    $\omega$& Maximal curvature & \\
    \hline
    $\varepsilon$ & Bending elastic constant & \multirow{3}{*}{Penalty parameter}\\
    $\nu$ & Curvature constraint elastic constant & \\
    $\mu $& Curvature control elastic constant & \\
    \hline
    $G$&$G[q,\nu,\varepsilon,\omega]=\varepsilon+\nu\left(|q_{ss}|^2-\omega^2\right)_+$
        & Reaction force\\
\hline
    $H$&$H[q,\mu,u,\omega]=\mu\left(\omega u-q_s \times q_{ss}\right)$    & Control term\\
        \hline
        $\beta$&Enviromental friction& \multirow{2}{*}{Friction coefficient}\\
        $\gamma$&Internal friction& \\
        \hline
        \end{tabular}}
   
\end{table} 
\begin{remark}[Deactivated controls]\label{parameterscontrol}
We recall that the control term in \eqref{sysintro} is the elastic force
$$F[q,\mu,u,\omega]:= - (H q_{ss}^\bot)_s -(H q_s^\bot)_{ss}$$
where $H=\mu(\omega u-q_s\times q_{ss})$. Therefore, to neglect the penalty parameter $\mu$ can  be used to model the deactivation of the controls in a subregion of the device. The scenarios we have in mind include mechanical breakdowns, voluntary deactivation for design or energy saving purposes and, as we see below, modeling applications to hyper-redundant systems. For instance, to choose $\mu\equiv 0$  yields the uncontrolled dynamics: 
\begin{equation}\label{uncontrolled}
\rho q_{tt}=(\sigma q_s)_s - (G q_{ss})_{ss}.\end{equation}
More generally, to set $\mu\equiv0$ in a closed set $I\subset[0,1]$ means that the portion of $q$ corresponding to $I$ is uncontrolled and it evolves according to \eqref{uncontrolled} -- with suitable (time dependent, controlled) boundary conditions.   
\end{remark}

We conclude this section by focusing on the tuning of the mass distribution $\rho$ and on the curvature constraint $\omega$ in order to encompass some morphological properties of the original three dimensional model.  
 Following \cite{CLL19}, we consider a three dimensional manipulator, endowed with axial symmetry and uniform mass density $\rho_v$. In particular, axial symmetry implies that the cross section of the manipulator, at any point $s\in[0,1]$  of its symmetry axis, is a circle $\Omega(s)$ of radius, say, $d(s)$. Therefore, to set $\rho(s):=\pi \rho_v d^2(s)$ corresponds to concentrate the mass of $\Omega(s)$ on its barycenter. 
The curvature constraint $\omega(s)$ can be chosen starting from the general consideration that a bending induces a deformation of the elastic material composing the body of the manipulator, for which assume an uniform yield point. In other words, in order to prevent inelastic deformations, the pointwise elastic forces acting on the material must be bounded by an uniform constant $F_{max}$. In view of the axial symmetry, the maximal angular elastic force $F(s)$ in any cross-section $\Omega(s)$ is attained on its boundary and
it reads $F(s)=e|q_{ss}|d(s)$, where $e$ is the elastic constant of the material.
Therefore to impose $F(s)\leq F_{max}$ is equivalent to the curvature constraint $|q_{ss}|\leq \omega(s):=F_{max}/(e d(s))$.

 \section{\uppercase{Optimal reachability}}\label{control}
 In this section we focus on an optimal reachability problem in both stationary and dynamic settings. The problem is to steer the end-effector (or tip) of the device, parametrized by $q(1,t)$, towards a target in the plane $q^*$, while optimizing accuracy, steadiness and energy consumption. Our approach is based on optimal control theory, that is, we recast the problem as a constrained minimization of a cost functional.
   \subsection{Static  optimal reachability}\label{stat}
   We begin by focusing on the optimal shape of the soft manipulator at the equilibrium. Equilibria $(q,\sigma)$ of the soft manipulator equation \eqref{sysintro} where explicitly characterized in \cite{CLL18}, and generalized to the case of uncontrolled regions of the manipulator in \cite{CLL19}. In particular, assuming the technical condition $\mu(1)=\mu_s(1)=0$, the shape of the manipulator  is the solution  $q$ of
the following second order ODE:
\begin{equation}\label{reducedstationary}
\quad\left\{\begin{array}{ll}
q_{ss}=\bar \omega u q_s^\bot&\mbox{in }(0,1)\\
|q_s|^2=1 &\mbox{in }(0,1)\\
q(0)=(0,0)\\q_{s}(0)=(0,-1).
\end{array}
\right. 	\end{equation} 
where  $\bar \omega:=\mu \omega/(\mu+\varepsilon)$.
\begin{remark}[Equilibria with deactivated controls]
 If $q$ is uncontrolled in $I\subset[0,1]$, i.e., $\mu(s)=0$ for all $s\in I$, then \eqref{reducedstationary} implies $|q_{ss}|=0$ in $I$, that is the corresponding portion of the device at the equilibrium is arranged in a straight line.
\end{remark}

We consider a cost functional involving a quadratic cost on the controls, modeling the energy required to force a prescribed curvature, and a tip-target distance term:
\begin{equation}\label{js}
\mathcal J^s:=\frac12\int_{[0,1]\setminus I} u^2 ds+\frac{1}{2\tau}|q(1)-q^\ast|^2,\end{equation}
where $I:=\{s\in(0,1)\mid \mu(s)=0\}$. Note that the domain of the control $u$ is in general the whole interval $[0,1]$ and thus it is independent from the parameter $\mu$, which is deputed to quantify the (possibly null) dynamic effects of $u$. To cope with this, the domain of integration of $\mathcal J^s$ is restricted to the regions in which the controls are really actuated. Note that, due to the quadratic dependence of the cost in $u$, we may equivalently extend the domain of integration of $\mathcal J^s$ to the whole $[0,1]$ by setting by default $u(s)\equiv 0$ for $s\in I$. Also remark that, tuning the penalization parameter $\tau$ allows to prioritize either the reachability task or the energy saving. 
The stationary optimal control problem then reads
\begin{equation}\label{touchfunctionalstationary}
\min \mathcal J^s,\quad\text{subject to \eqref{reducedstationary} and to $|u|\leq 1$.}
\end{equation}
Following \cite{CLL19}, we can restate \eqref{touchfunctionalstationary} in terms of the Euler's elastica type variational problem
\begin{equation}\label{staticfunctional}
 \begin{split}\min\left\{\frac12\int_{(0,1)\setminus I} \frac{1}{\bar\omega^2}|q_{ss}|^2 ds+\frac{1}{2\tau}|q(1)-q^\ast|^2\right\}
 \\
 \mbox{ subject to }\left\{\begin{array}{ll}
 |q_s|^2=1 &\mbox{in }(0,1)\\
 |q_{ss}|\leq \bar \omega&\mbox{in }(0,1)\\
 q(0)=(0,0)\\q_{s}(0)=(0,-1).
 \end{array}
 \right.\end{split}
 \end{equation}
 Indeed, by differentiating the constraint $|q_s|=1$ we obtain the relation $q_{ss}\cdot q_s=0$ and, consequently, $|q_{ss}\cdot q_s^\bot|=|q_{ss}|\cdot |q_s^\bot|=|q_{ss}|$. Then, by dot multiplying $q^\bot$ in both sides of the first equation of \eqref{reducedstationary}, we deduce $|q_{ss}|=\bar\omega |u|$ and, consequently, 
 $u^2=\frac{1}{\bar \omega^2}|q_{ss}|^2$ in $(0,1)\setminus I$.

 \subsubsection{Numerical tests for static optimal reachability.} 
%  We complete this section by remarking how tuning the parameter $\mu$ may used to model the impact of the actuation of the controls on the system. In particular in some regions $I\subset[0,1]$ we may not want or can control the device: this can be modelled by setting the control penalty parameter $\mu:=\mu_I\equiv0$ in $I$.
 We consider numerical solutions of the optimal problem \eqref{staticfunctional} in the following scenarios:  $I=\emptyset$, i.e., the deviced is fully controlled, $I=[0.35,0.65]$, i.e., the device is uncontrolled in its median section,  $I=[0.25,0.4]\cup[0.6,0.75]$ and, finally, $I=[0,1]\setminus\{0,0.25,0.5,0.75\}$. Slight variations (in terms of dicretization step) of the first three tests were earlier discussed in \cite{CLL19},  while the last test is completely new.  
 
       %To numerically solve \eqref{staticfunctional} we encompass the constraints $|q_s|=1$ and $|q_{ss}|\leq \bar \omega$ %(encompassing the fact that the control set for $u$ is the bounded interval $[-1,1]$)  via an augmented Lagrangian method, and we discretize it with a finite difference scheme. The non-linear terms of the resulting discrete system of equations are treated via a quasi Newton's method.  Parameters settings are reported in Table \ref{parametertable}.
   \begin{table}[t]	\caption{Control deactivation settings related to $\mu_I$. \label{controltable}}\centering
  {	\renewcommand\arraystretch{1.3}\begin{tabular}{|c|l|}\hline
  		{\bf Test}& {\,\,\bf Control deactivation region\,\,}\\
  		\hline
  	\,\,	Test 1\,\,& \,\,$I=\emptyset$\\
  	\,\,	Test 2\,\,& \,\,$I=[0.35,0.65]$\\
  	\,\,	Test 3\,\,& \,\,$I=[0.25,0.4]\cup[0.6,0.75]$\\
  	\,\,	Test 4\,\,& \,\,$I=[0,1]\setminus\{0,0.25,0.5,0.75\}$\\
  	\hline
  	\end{tabular}}
  
  \end{table}
   \begin{table}[ht]\caption{Global parameter settings.  \label{parametertable}}\centering{\renewcommand{\arraystretch}{1.3} 
  	\begin{tabular}{|l|l|}\hline
  {\bf\,\,Parameter description}\,\,&{\bf\,\,Setting}\\
  		\hline
  		\,\,Bending moment &\,\,$\varepsilon(s)=10^{-1}(1-0.9s) $\\
  		  	\,\,Curvature control penalty&\,\,$\mu(s)=(1-s) \exp(-0.1\frac{s^2}{1-s^2})$\,\,\\
  	  		\,\,Curvature constraint&\,\,$\omega(s)=4\pi(1+s^2)$\\
  		\hline
  		\,\,Target point&  $q^*=(0.3563,-0.4423)$\\
  	\,\,Target penalty& $\tau=10^{-4}$\\
  		\hline
  		\,\,Discretization step& $\Delta_s=0.02$\\\hline
  	\end{tabular}}
  	  \end{table}
 
%  \begin{figure}[ht]
%   	\centering{ \scalebox{0.9}{	\begin{tabular}{cc}
%   		\includegraphics[width=0.43\textwidth]{shape-01}& \includegraphics[width=0.523\textwidth]{curvature-01}\\
%   		(a)&(b)	
%   	\end{tabular}}
%   	\caption{\label{stationary1} In (a) the solution $q$ of Test 1, in (b) the related signed curvature $\kappa(s)$ (bold line) and curvature constraints $\pm\bar \omega$ (thin lines).}}
%   \end{figure} 
   \begin{figure}[h!]
	\centering \scalebox{1}{\renewcommand{\arraystretch}{1.3}	
	\begin{tabular}{ |c|c|}
	\hline$q(s)$&$\kappa(s)$\\\hline
	 \includegraphics[width=.43\textwidth]{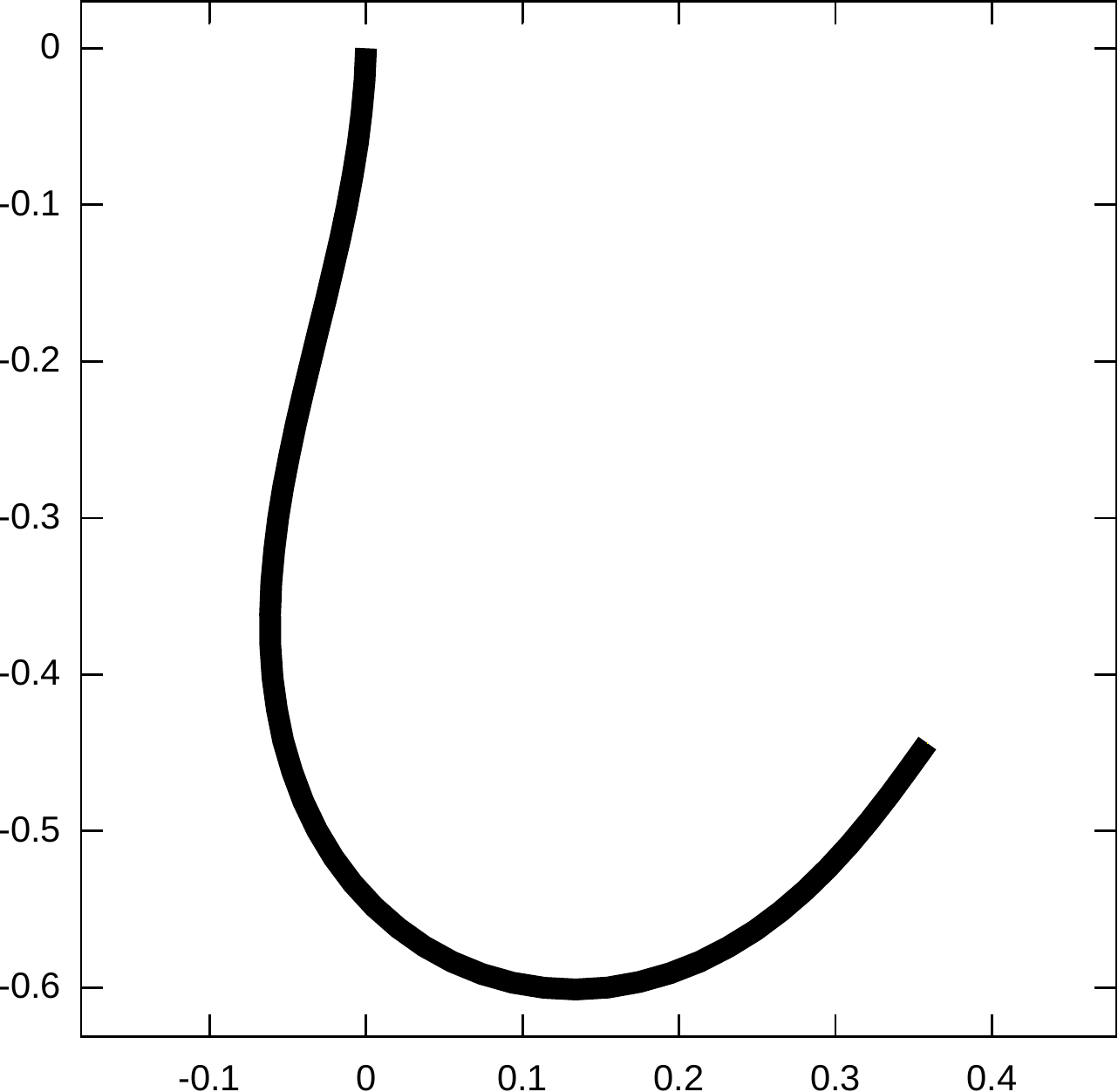}\,\,\,&
		\includegraphics[width=.523\textwidth]{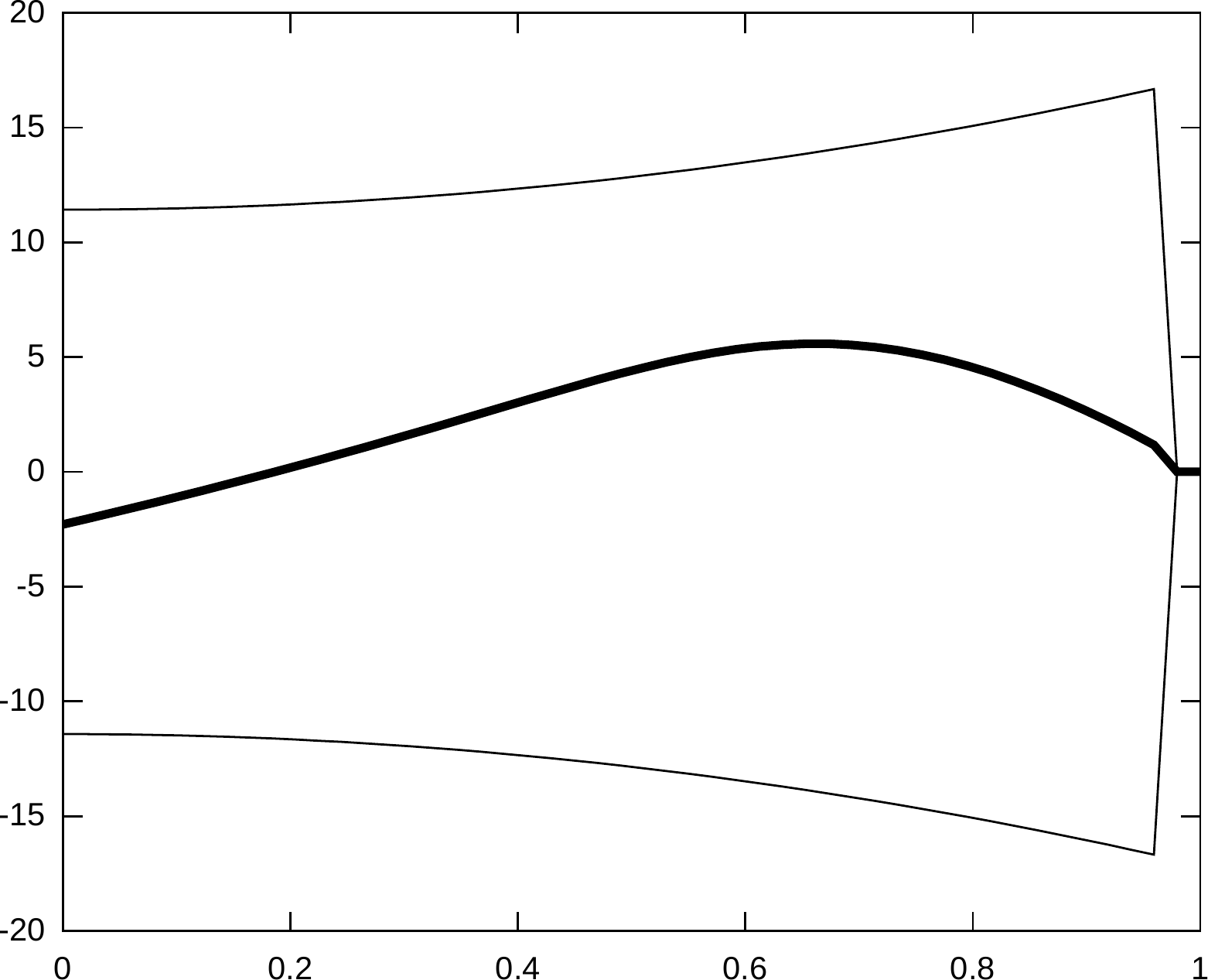}\,\,\,\\
% 	(a1)&(b1)\\
	 \includegraphics[width=.43\textwidth]{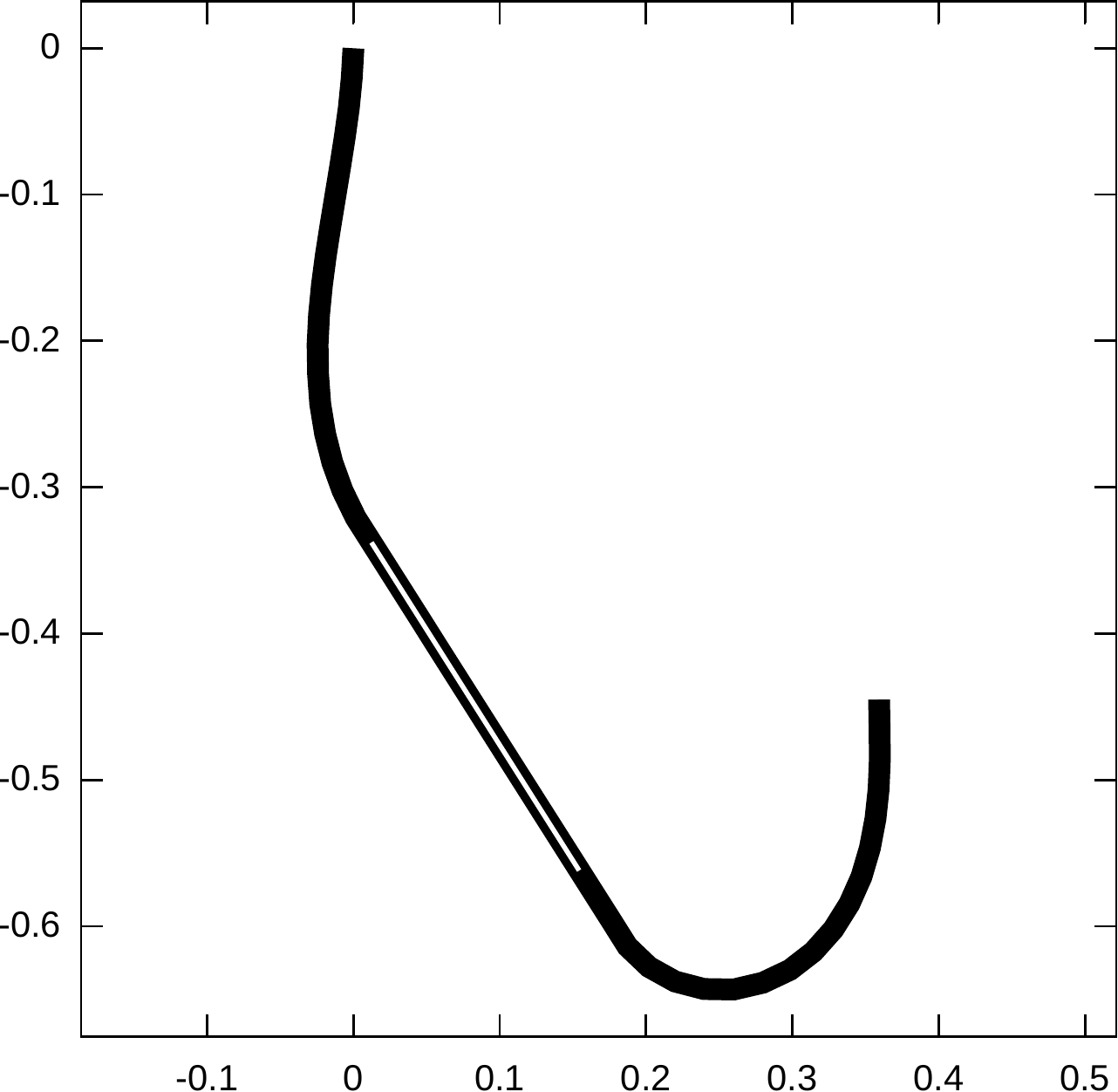}\,\,\,&
		\includegraphics[width=.523\textwidth]{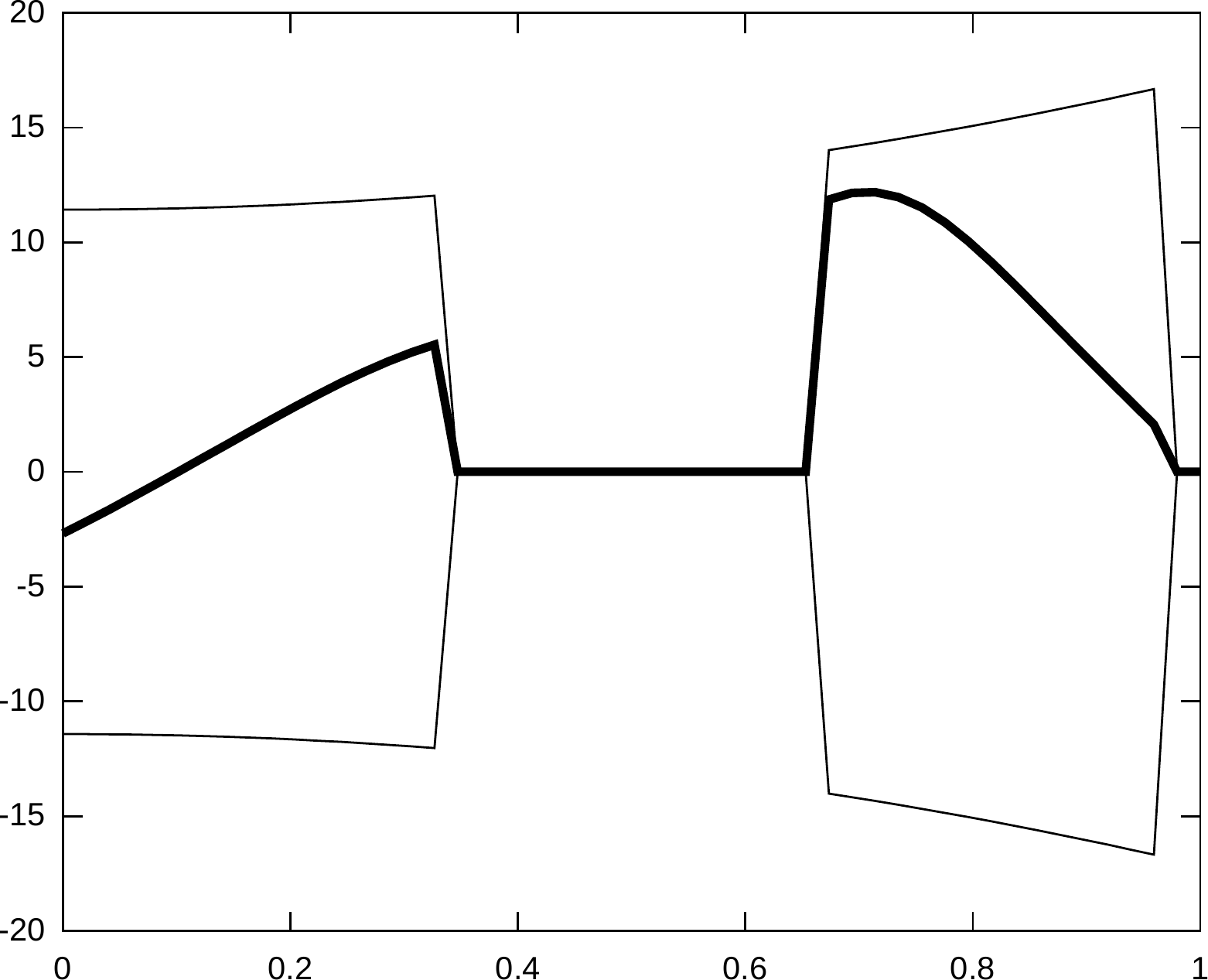}\,\,\,\\
% 			(a2)&(b2)	\\
	\includegraphics[width=.43\textwidth]{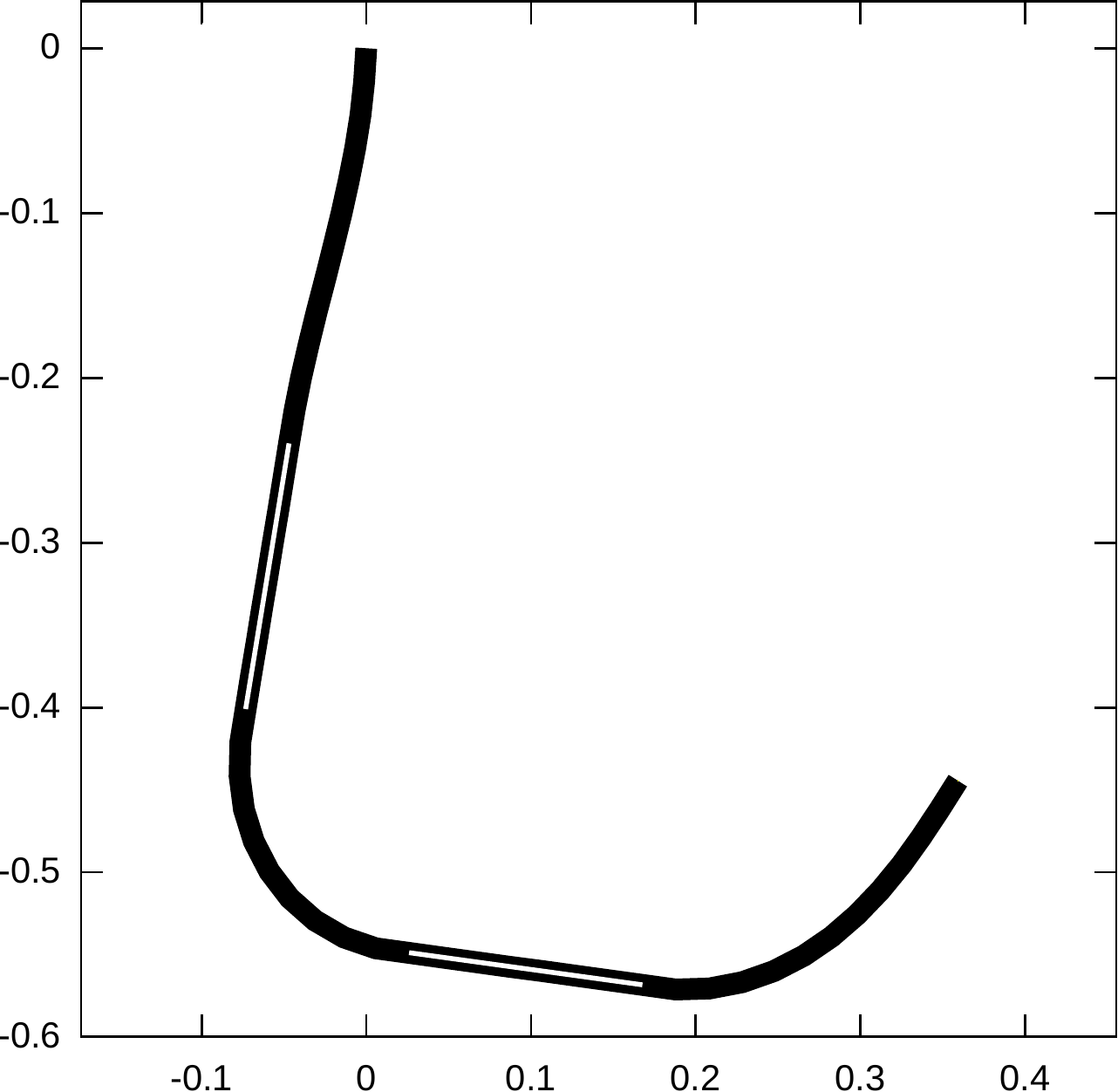}\,\,\,&
		\includegraphics[width=.523\textwidth]{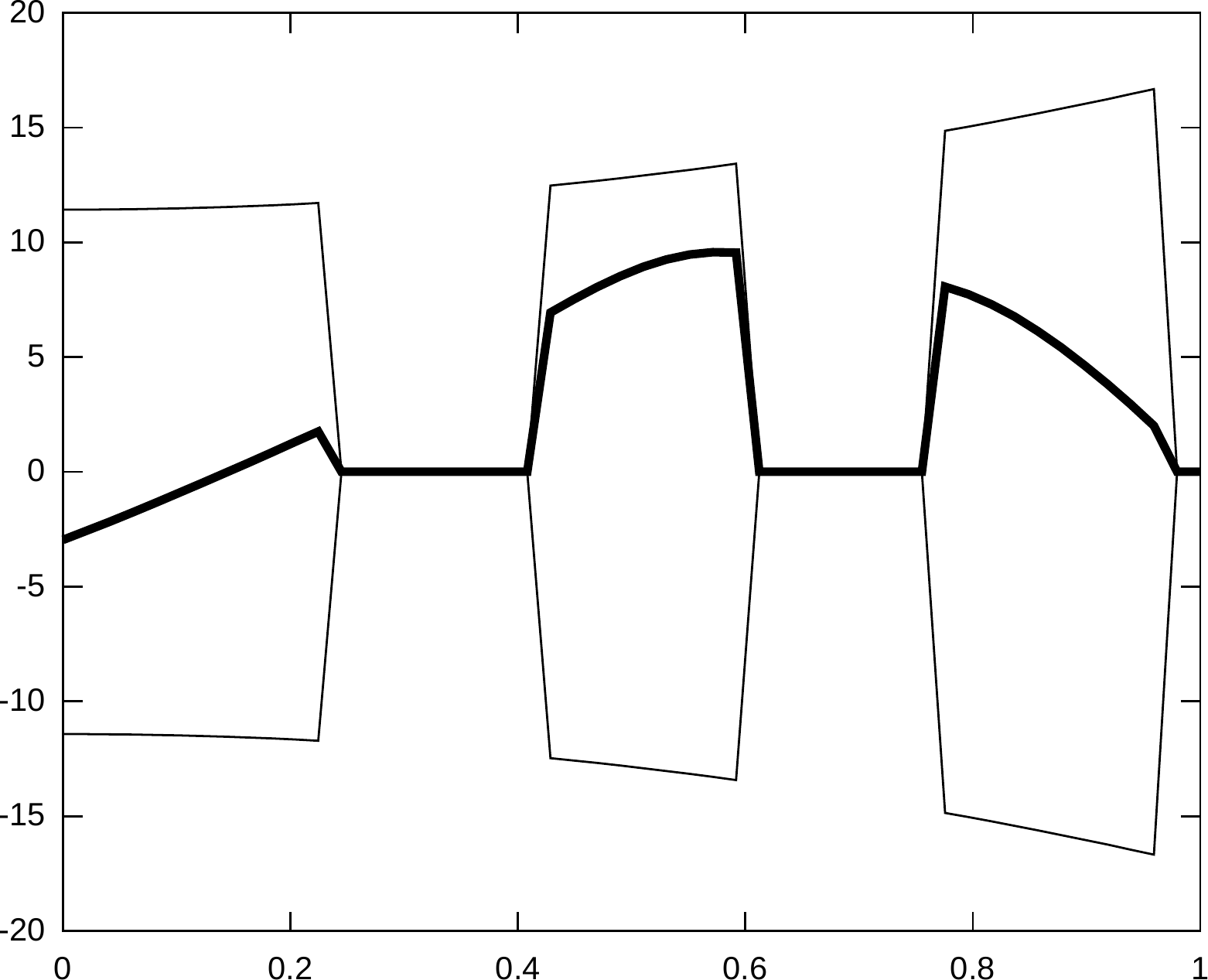}\,\,\,\\\hline
% 			(a3)&(b3)
			\end{tabular}}
	\caption{\label{stationary1}  On the first, second and third row are reported the results of Test 1, 2 and 3, respectively. In the first column are depicted the optimal solutions and, in the second column, the related signed curvatures (bold line) and curvature constraints $\pm\bar \omega$ (thin lines).} 
\end{figure}

  \begin{figure}[h!]
	\centering \scalebox{0.99}{	
	\begin{tabular}{cc}
	{\centering$q(s)$}&$\kappa(s)$\\
 \includegraphics[width=.43\textwidth]{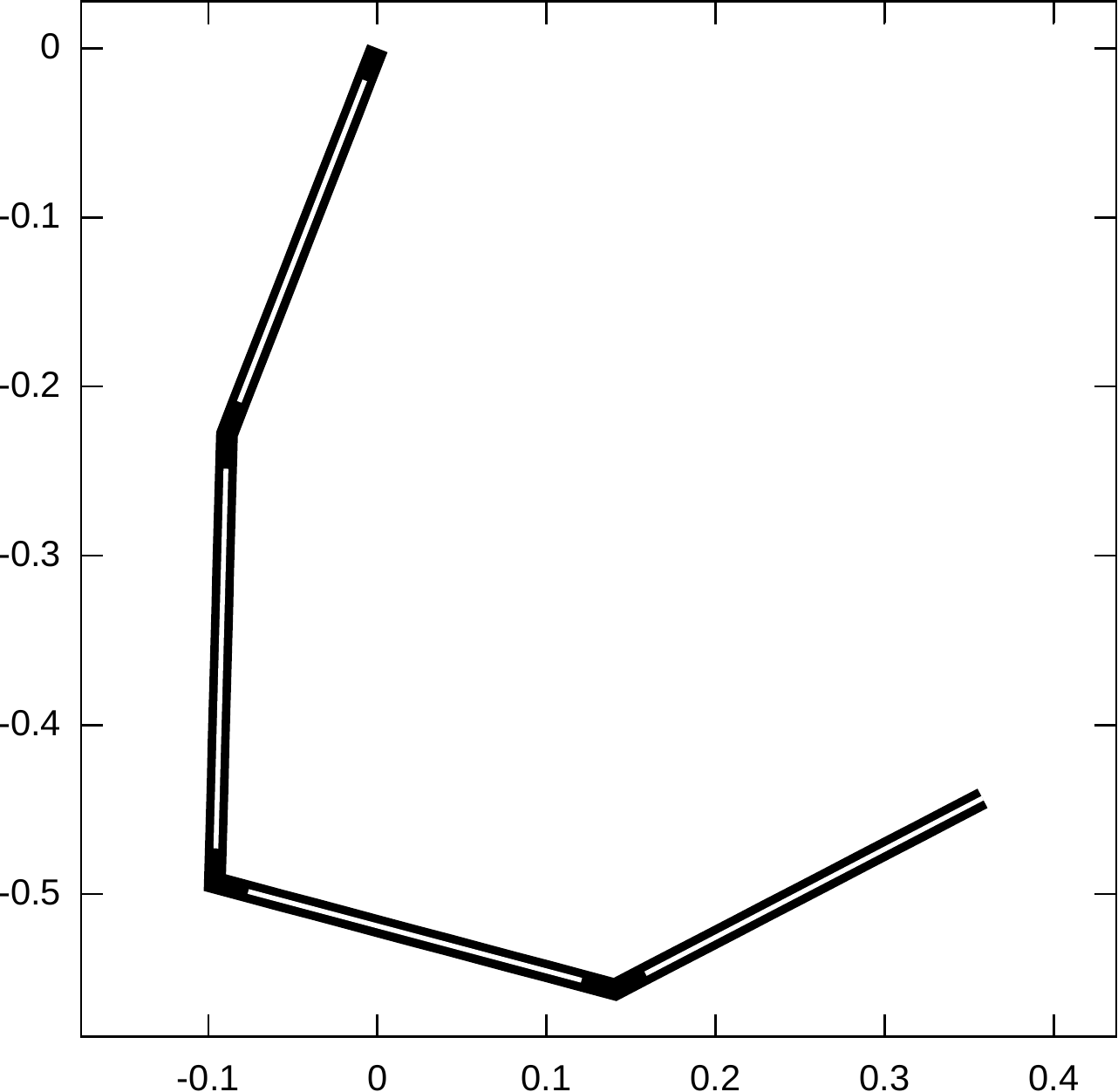}&
		\includegraphics[width=.523\textwidth]{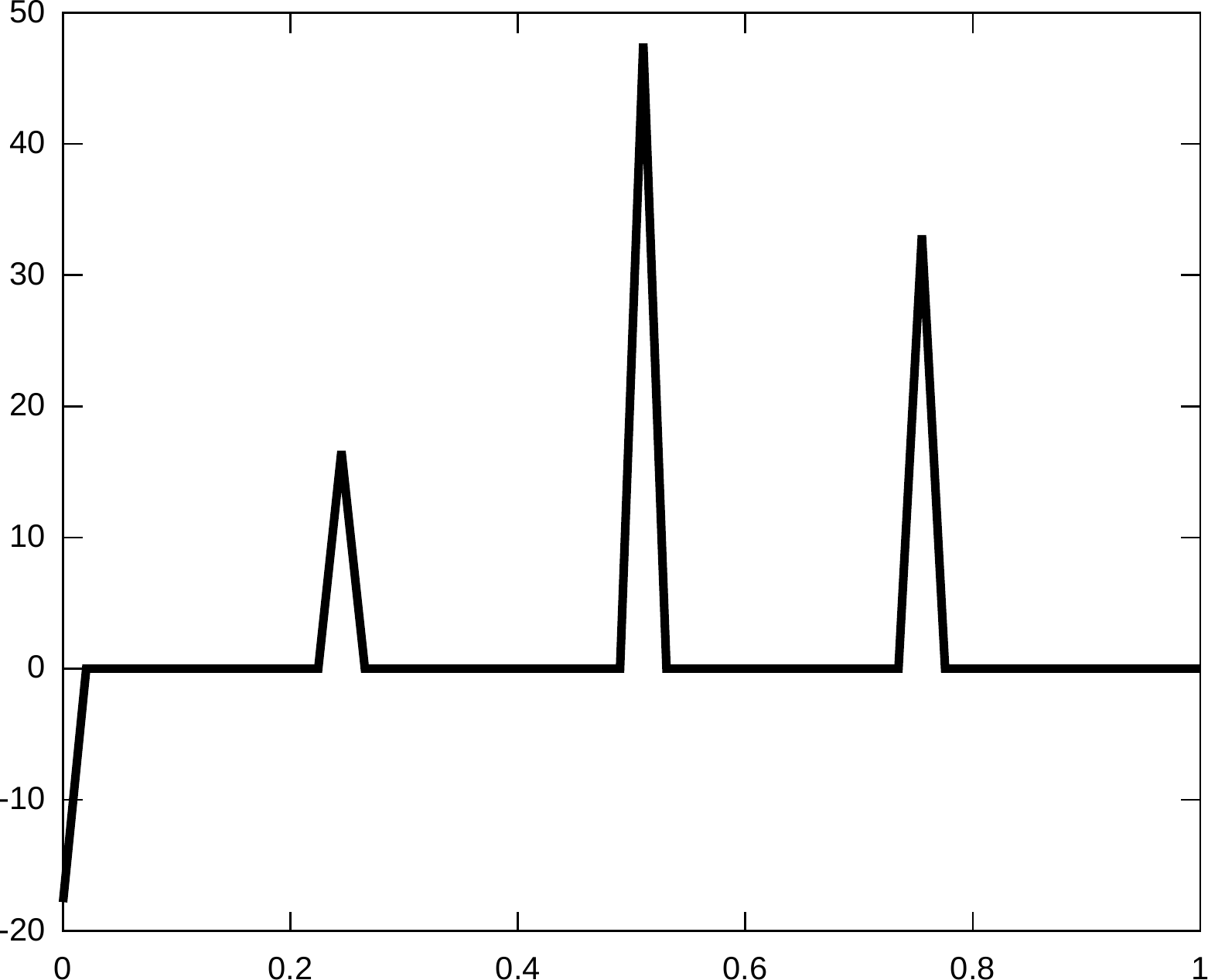}\\(a)&(b)	\end{tabular}}
	\caption{\label{stationary4} In (a) the solution $q(s)$ of Tests 4 and in (b) the related signed curvature $\kappa(s)$.} 
\end{figure}

In the numerical implementation, the constraints $|q_s|=1$ and $|q_{ss}|\leq \bar \omega$ %(encompassing the fact that the control set for $u$ is the bounded interval $[-1,1]$) 
are enforced via an augmented Lagrangian method, and the problem is discretized  using a finite difference scheme. This  yields a discrete system of equations, whose non-linear terms are treated via a quasi Newton's method. 
 
In what follows we depict  the solution $q(s)$ of the problem \eqref{staticfunctional} and the associated signed curvature $\kappa(s):=q_s(s)\times q_{ss}(s)$
in the cases reported in Table \ref{controltable}. Note that, due to the choice of arclenght coordinates, the (unsigned) curvature of $q$ is $|\kappa|=|q_{ss}|$.

We report the results of Test 1, Test 2 and Test 3 --see Table \ref{parametertable}-- in   Figure \ref{stationary1}: regions corresponding to in-actuated control are depicted in white. 
  
  Test 4 requires a few more words of discussion. In this case the manipulator is completely uncontrolled, with the exception of four points. The aim is to recover, from a merely static point of view, optimal configurations for a (hyper-redundant) rigid manipulator from the soft manipulator model: the controlled points  represent indeed the joints of the manipulator. Using the fact that uncontrolled regions arrange in straight lines at the equilibrium, we obtain an optimal configuration modeling at once the equilibrium of the soft manipulator as well as the solution of an optimal inverse kinematic problem for a hyper-redudant manipulator -- see Figure \ref{stationary4}-(a). Note that concentrating the actuation of controls on singletons yields the support of the resulting optimal curvature to concentrate on the same isolated points -- see Figure \ref{stationary4}-(b). To cope with this expected phenomenon, in Test 4 the curvature constraint is dropped.

 We finally note that the particular choice of $q^*$ guarantees reachability in all these four cases, but clearly, with different optimal solutions.

%   The discretization step is $\Delta=??$, i.e., the system is approximated as a system of $N=??$ particles. Note that the discretization  can also physically interpreted a system of $N-1$ rigid joints, subject to  (properly scaled) constraints on their relative rotations. The arrest criterion for the algorithm is given by the global tolerance $tol=??$ (i.e., the algorithm stops when two consecutive steps yield a variation in the discretized cost functional  less than $tol$) whereas Newton's step tolerance is $tol_{Newtown}=??$. 

 \subsection{Dynamic optimal reachability}
 We now consider the dynamic version of the reachability problem discussed in Section \ref{stat}: given $q^\ast\in\RR^2$ and $T, \tau>0$, we look for a \emph{time-varying} optimal control $u^*:[0,1]\times[0,T]\to[-1,1]$ minimizing the functional
 \begin{equation}\label{dynamicfunc}
 \begin{split}\mathcal{J}(q,u)&=\frac{1}{2\tau}\int_0^T\hskip-5pt|q(1,t)-q^\ast|^2 dt+\frac{1}{2}\int_0^T\hskip-5pt\int_0^1\hskip-5pt u^2 ds\,dt\\
 &+\frac{1}{2}\int_0^1\hskip-5pt\rho(s)|q_t(s,T)|^2 ds,
 \end{split}\end{equation}
 subject to the soft manipulator dynamics \eqref{sysintro}. The terms representing the tip-target distance and the quadratic cost on the control are now declined in a dynamic sense: they are minimized during the whole evolution of the system. The third term of $\mathcal J$ is deputed to steadiness, it represents the kinetic energy of the whole manipulator at final time. 
  Our approach is based on the study of first-order optimality conditions, i.e., on the solution of the so-called \emph{optimality} system. 
  The unknowns are the stationary points $(q,\sigma,u)$ of $\mathcal J$ (subject to \eqref{sysintro}), and the related multipliers $(\bar q,\bar \sigma)$ which are called the \emph{adjoint states}.  Roughly speaking, 
 if a control $u^*$ is optimal, then the optimality system admits a solution of the form $(q,\sigma,u^*,\bar q,\bar \sigma)$.  
  The optimality system is composed by two PDEs, describing the evolution of the adjoint states, the equations of motion \eqref{sysintro} and a variational inequality for the control.
  
   Following \cite{CLL18,CLL19},  the adjoint
  states equations are: 
  \begin{equation}\label{adjointfriction}
 \begin{cases}
 \rho \bar q_{tt}=\left(\sigma \bar q_s-H \bar q_{ss}^\bot\right)_s -\left(G \bar q_{ss}+H \bar q_{s}^\bot\right)_{ss} \\
 \qquad +\left(\bar\sigma  q_s-\bar H q_{ss}^\bot\right)_s -\left(\bar G q_{ss}+\bar H q_{s}^\bot\right)_{ss}\\
 \qquad +\beta \bar q_t+\gamma \bar q_{sssst}\\
 \bar q_s \cdot q_s=0 
 %\\\\
 %\rho q_{tt}=\left(\sigma q_s-H q_{ss}^\bot\right)_s -\left(G q_{ss}+H q_{s}^\bot\right)_{ss}\\
 %|q_s|^2=1   
 \end{cases}
 \end{equation}
 for $(s,t)\in(0,1)\times(0,T)$. The maps $G$ and $H$ are defined in Section \ref{model}
 and the maps $\bar G$ and $\bar H$ are their linearizations, respectively:
 	$$\bar G[q,\bar q,\nu,\omega]=g[q,\nu,\omega]q_{ss}\cdot \bar q_{ss}\,,$$
 		$$\bar H[q,\bar q,\mu]=\mu \left(\bar q_s \times q_{ss}+q_s \times \bar q_{ss}\right)\,,$$
 		 	where $g[q,\nu,\omega]=2\nu\mathbf{1}(|q_{ss}|^2-\omega^2)$ and $\mathbf{1}(\cdot)$ stands for the Heaviside function, i.e. $\mathbf{1}(x)=1$ for $x\ge 0$ and $\mathbf{1}(x)=0$ otherwise. 	 
 		 	
 Optimality conditions also yield final and boundary conditions on $(\bar q,\bar \sigma)$. As a consequence of the fact that the optimization takes into account the whole time interval $(0,T)$, initial conditions on $q$ correspond to final conditions on its adjoint state $\bar q$:
 $$\bar q(s,T)=-q_t(s,T), \qquad\bar q_t(s,T)=0\qquad \text{for }s\in (0,1).$$	
 Boundary conditions are reported in Table \ref{tboundary}. Note that the  zero, the first and the second order conditions on $\bar q$ display an essential symmetry with those on $q$ reported in \eqref{bond1}. On the other hand, both the third order condition on $\bar q$  and the adjoint tension boundary condition on $\bar \sigma$ show a dependence on the difference vector between the tip and the target $q^*$: this phenomenon represents the fact that the tip is forced towards $q^*$.
 \begin{table}[t]
 \caption{\label{tboundary}Boundary conditions on $(\bar q,\bar \sigma)$ for $t\in(0,T)$. 
 }
 \centering{\renewcommand{\arraystretch}{1.3}
 \begin{tabular}{|l|l|}\hline
 \multirow{2}{*}{\,\,Fixed endpoint} \,\,& \,\,$\bar q(0,t)=0$     \\
  &$\bar q_s(0,t)=0$    \\
  \hline
  \multirow{3}{*}{\,\,Free endpoint\,\,}&\,\,$\bar q_{ss}(1,t)=0$ \\
  &\,\,$\bar q_{sss}(1,t)=\frac{1}{\tau\varepsilon}\left((q-q^*)\cdot q_s^\bot\right) q_s^\bot (1,t)$\,\,\\
  &\,\,$ \bar\sigma(1,t)=-\frac{1}{\tau}(q-q^\ast)\cdot q_s(1,t)$\,\,
  \\\hline
 \end{tabular}}
  \end{table}
 
Finally, the variational inequality for the control is 
\begin{equation}\label{varinequality}
\int_0^{T}\int_0^1 \left( u + \omega\bar H[q,\bar q]\right) (v-u)ds\,dt\ge 0\,.
\end{equation}
for every $v:[0,1]\times[0,T]\to[-1,1]$, which provides, in a weak sense, the variation of the
functional $\mathcal J$ with respect to $u$, subject to the constraint $|u|\leq 1$.
% \section{\uppercase{Numerical simulations}}\label{numerics}

We remark that, due to friction forces, the system \eqref{sysintro} is dissipative. This implies that if we plug  in \eqref{sysintro} a constant in time control  $u^s$, that we call \emph{static optimal control}, given by the solution of the stationary control problem \eqref{touchfunctionalstationary}, then the system converges as $T\to +\infty$ to the optimal stationary equilibrium $q$ described in Section \ref{stat} - see  \eqref{reducedstationary}. 
\subsubsection{Numerical tests for dynamic optimal reachability.}In the following simulations optimal static controls are used as initial guess for the search of dynamic optimal controls, as well as benchmarks for dynamic optimization. In particular, the optimality system is discretized using a finite difference scheme in space, and a velocity Verlet scheme in time. Then the optimal solution is computed iteratively, using an adjoint-based gradient descent method. More precisely, starting from the optimal static control $u^s$ as initial guess for the controls, we first solve \eqref{sysintro} forward in time. The solution-control triplet $(q,\sigma,u)$ is plugged into \eqref{adjointfriction}, which in turn is solved backward in time. We obtain then the solution-control vector $(q,\sigma,u,\bar q,\bar \sigma)$ which is used in \eqref{varinequality} to update the value of $u$. This routine is iterated up to convergence in $u$. 

%-- incidentally remark that $s\in I$ implies $\bar\omega(s)=0$ which in turn implies $|q_{ss}|=0$ and, consequently, an indeterminate form of type $\frac{0}{0}$ in the integral in \eqref{staticfunctional}. 
    \begin{table}[t]\caption{Dynamic parameter settings.  \label{dynamictable}}\centering{\renewcommand{\arraystretch}{1.3}\begin{tabular}{|l|l|}
 		\hline
 		{\bf\,\,Parameter description}\,\,&\,\,{\bf Setting}\\\hline
 		\,\,Mass distribution&\,\,$\rho(s)=\exp(-s)$\\
 		  		\,\,Curvature constraint penalty \,\,&\,\,$\nu(s)=10^{-3}(1-0.09s)$\\
 		\,\,Environmental friction &$\,\,\beta(s):=2-s$\\
 		\,\,Internal friction &\,\,$\gamma(s):=10^{-6}(2-s)$\,\,\\
 		\,\,Final time &\,\,$T=2$\\
 		\,\,Time discretization step &\,\,$\Delta_t=0.001$\\\hline
 	\end{tabular}}
 \end{table}

%  The adjoint system composed by \eqref{sysintro}, \eqref{adjointfriction} and \eqref{varinequality} can be discretized using a standard finite difference scheme in space-time, then solved by an adjoint-based gradient descent method. The key idea is the following: starting from an initial guess $u$ given by the stationary optimal control, we first solve the equation of motion \eqref{tentaclemotionfriction} forward in time. Then the solution-control triplet $(q,\sigma,u)$ is plugged into \eqref{adjointfriction} which is  solved backward in time. Finally, we use the vector $(q,\sigma,u,\bar q,\bar \sigma)$ in \eqref{varinequality} to update the value of the $u$.  The procedure is iterated up to convergence on the control. 

We recall from \cite{CLL19} the investigation of the dynamic counterpart of Test 2,  a scenario in which the median section of the device is uncontrolled -- see Section \ref{stat} and, in particular, Table \ref{controltable}. Parameters are set in Table \ref{parametertable} and, for the dynamic aspects, in Table \ref{dynamictable}. In what follows, we compare the performances of static and dynamic optimal controls. In particular, we consider the dynamic optimal control $u^d$, i.e., the numerical solution of \eqref{dynamicfunc}, and the static optimal control $u^s$, which we recall is constant in time and it coincides with the solution of \eqref{staticfunctional}, see Figure \ref{stationary1}. 
 \begin{figure}[t]
   	\centering
   	\begin{tabular}{|c|c|c|}\hline
   	$\mathcal J_{q^*}$&$\mathcal J_{u}$&$\mathcal J_{v}$\\\hline
   		\includegraphics[width=.32\textwidth]{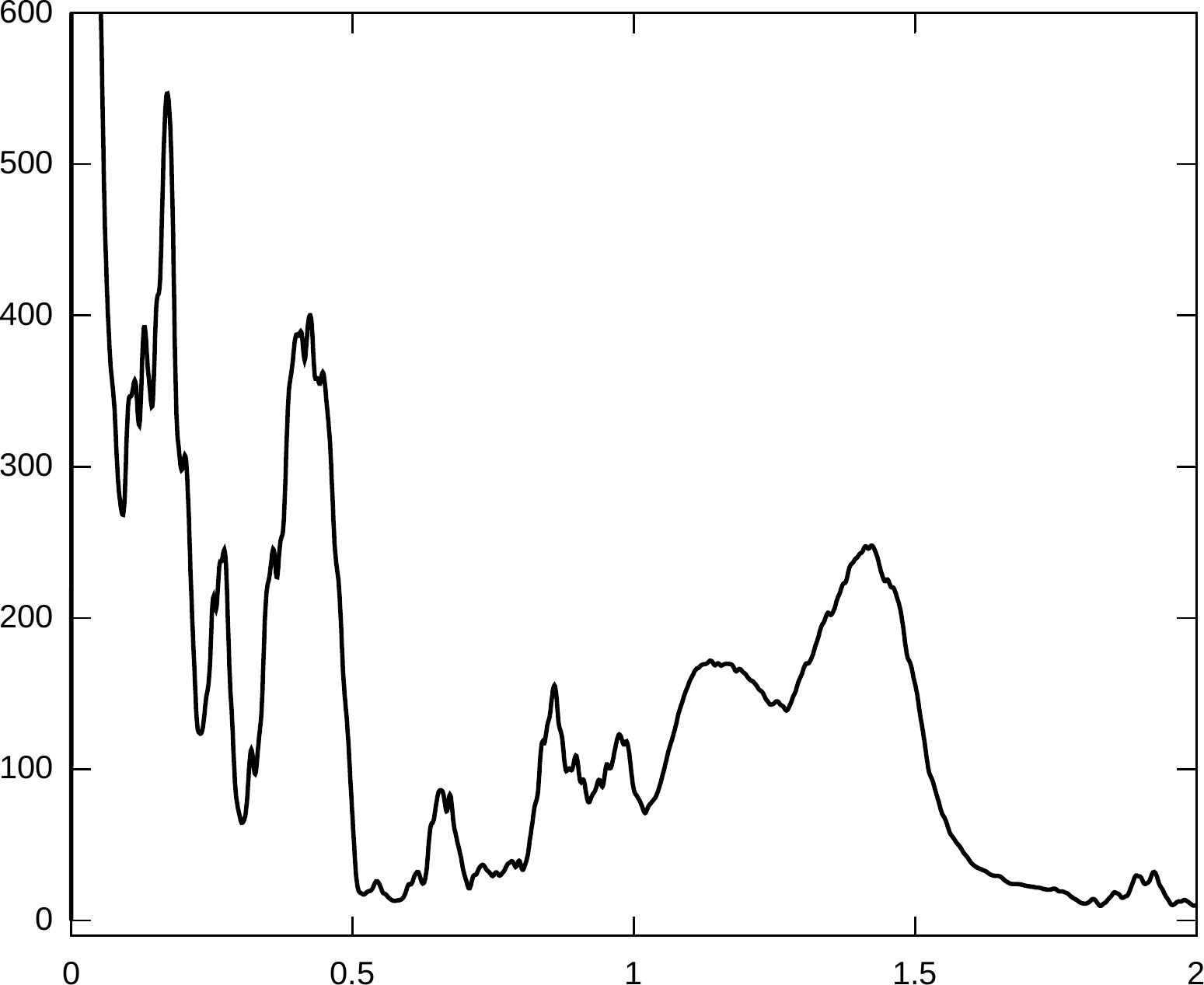}&\includegraphics[width=.32\textwidth]{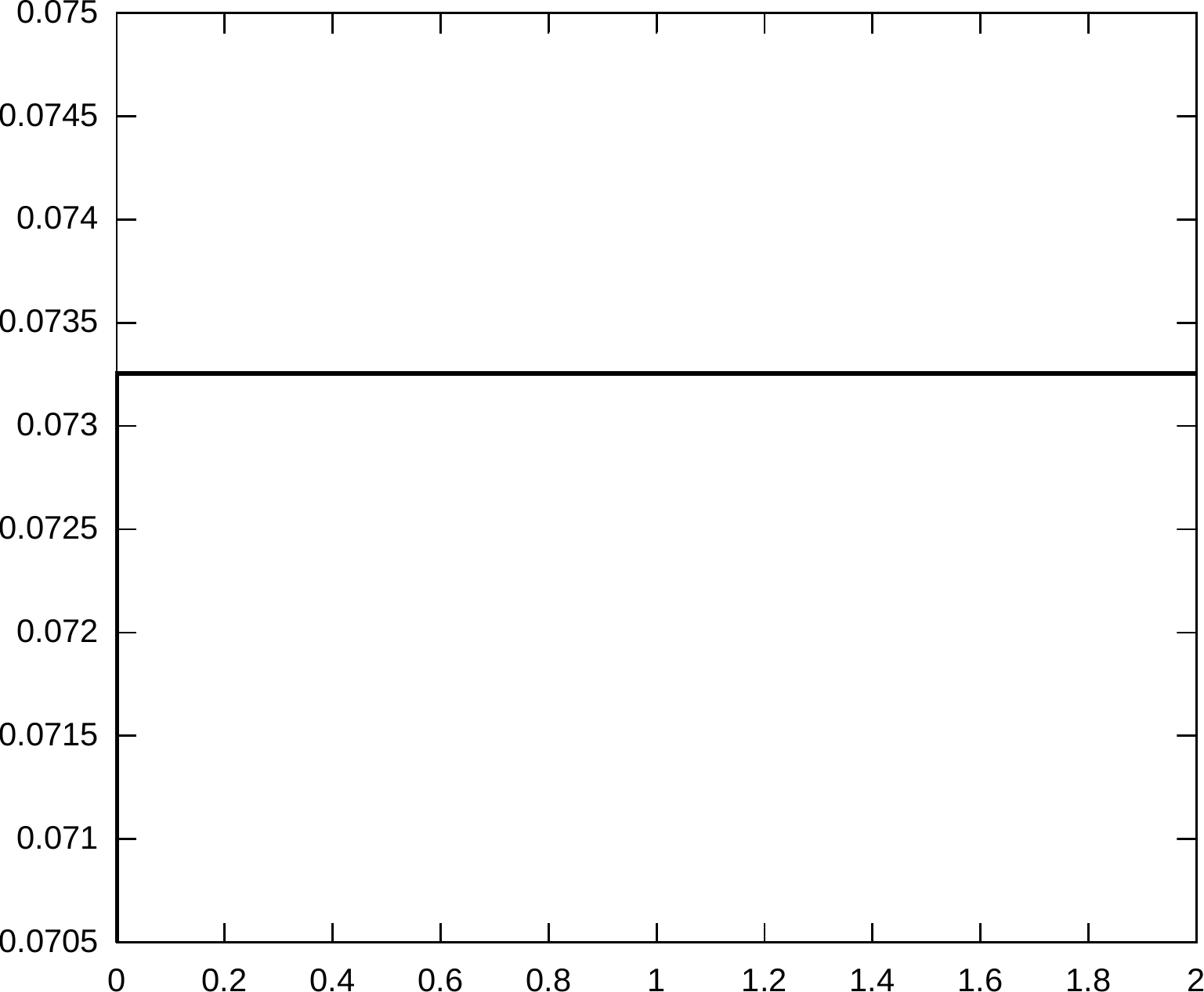}&\includegraphics[width=.32\textwidth]{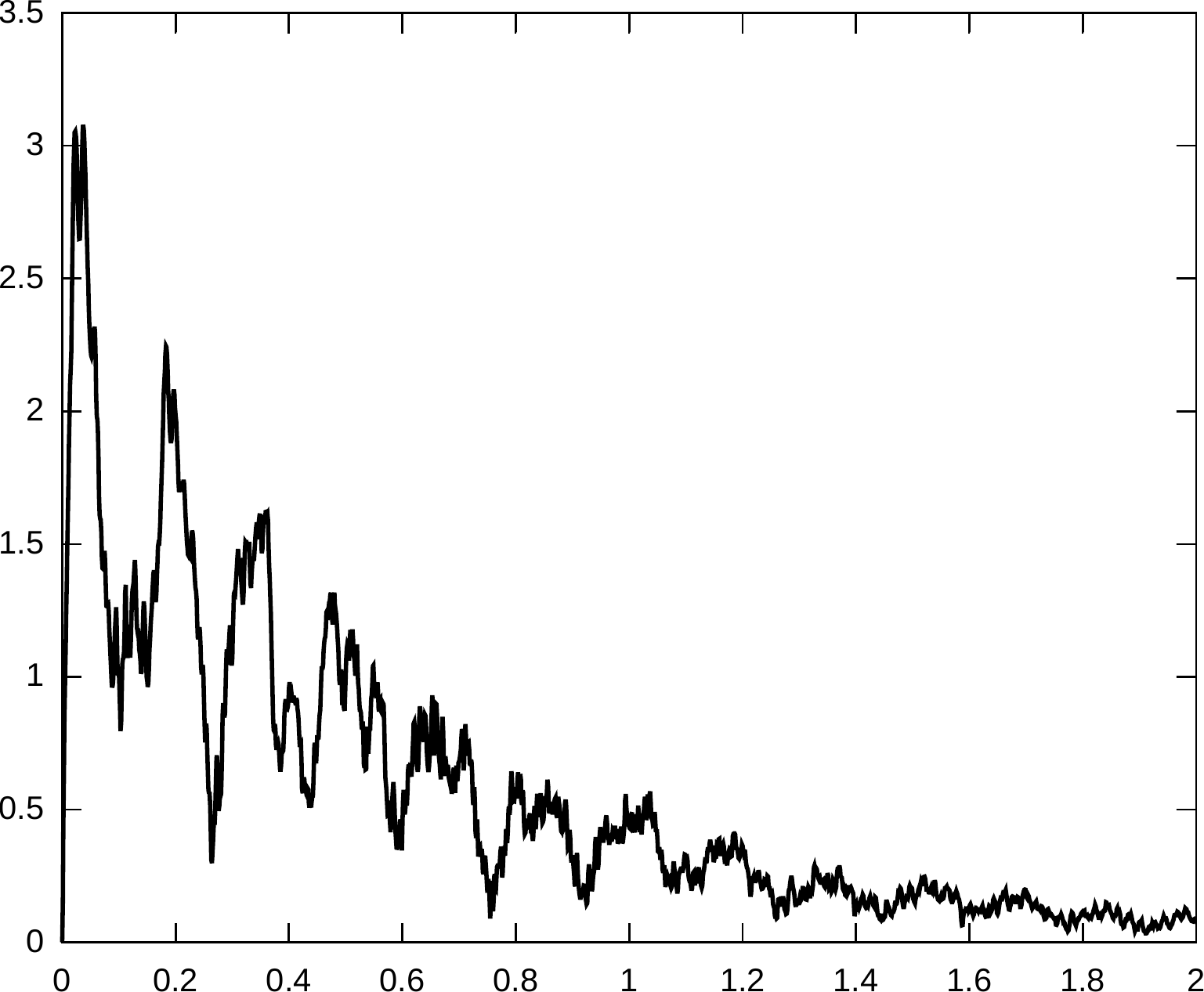}\\
   		&&\\
	\includegraphics[width=.32\textwidth]{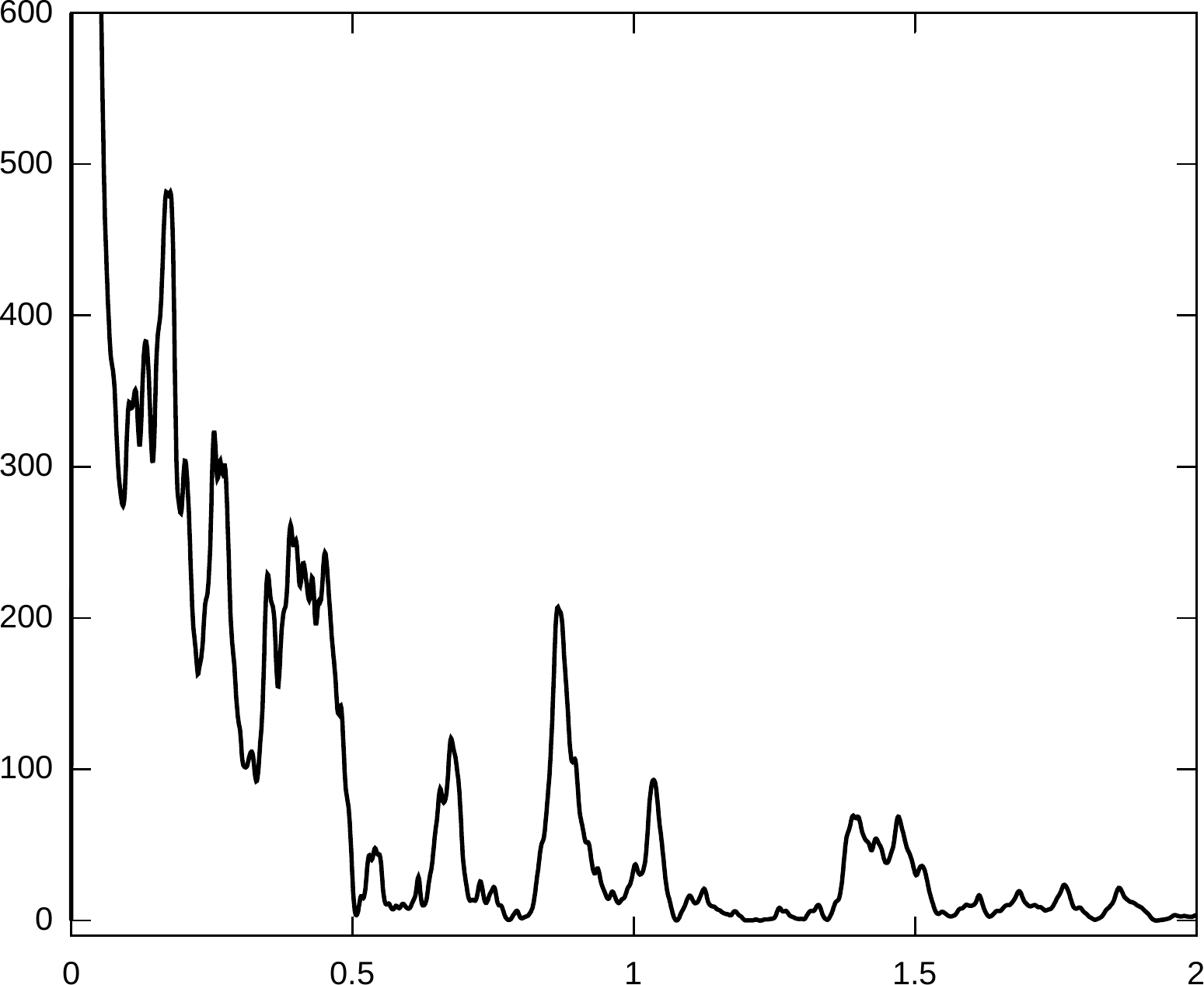}&\includegraphics[width=.32\textwidth]{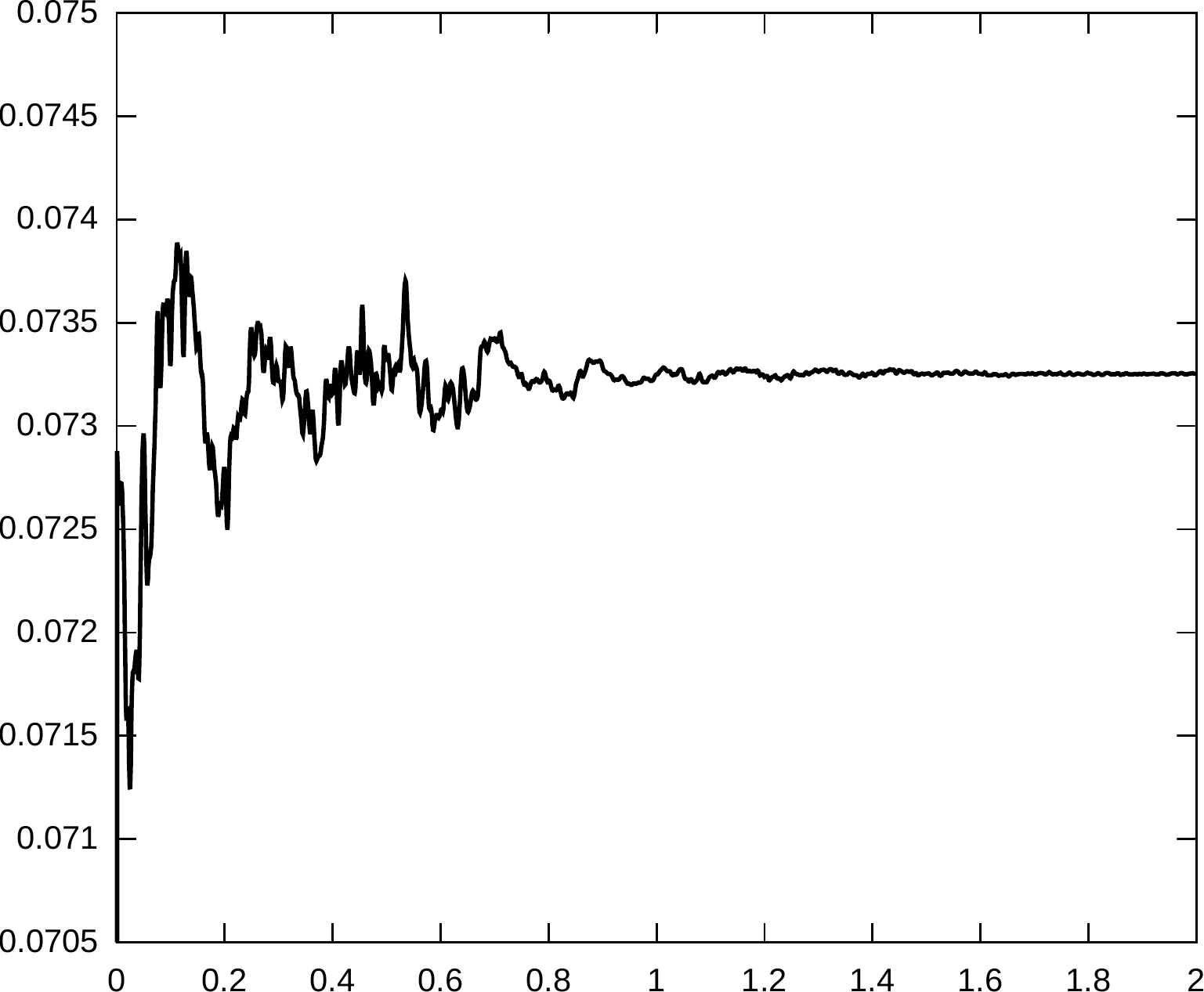}&\includegraphics[width=.32\textwidth]{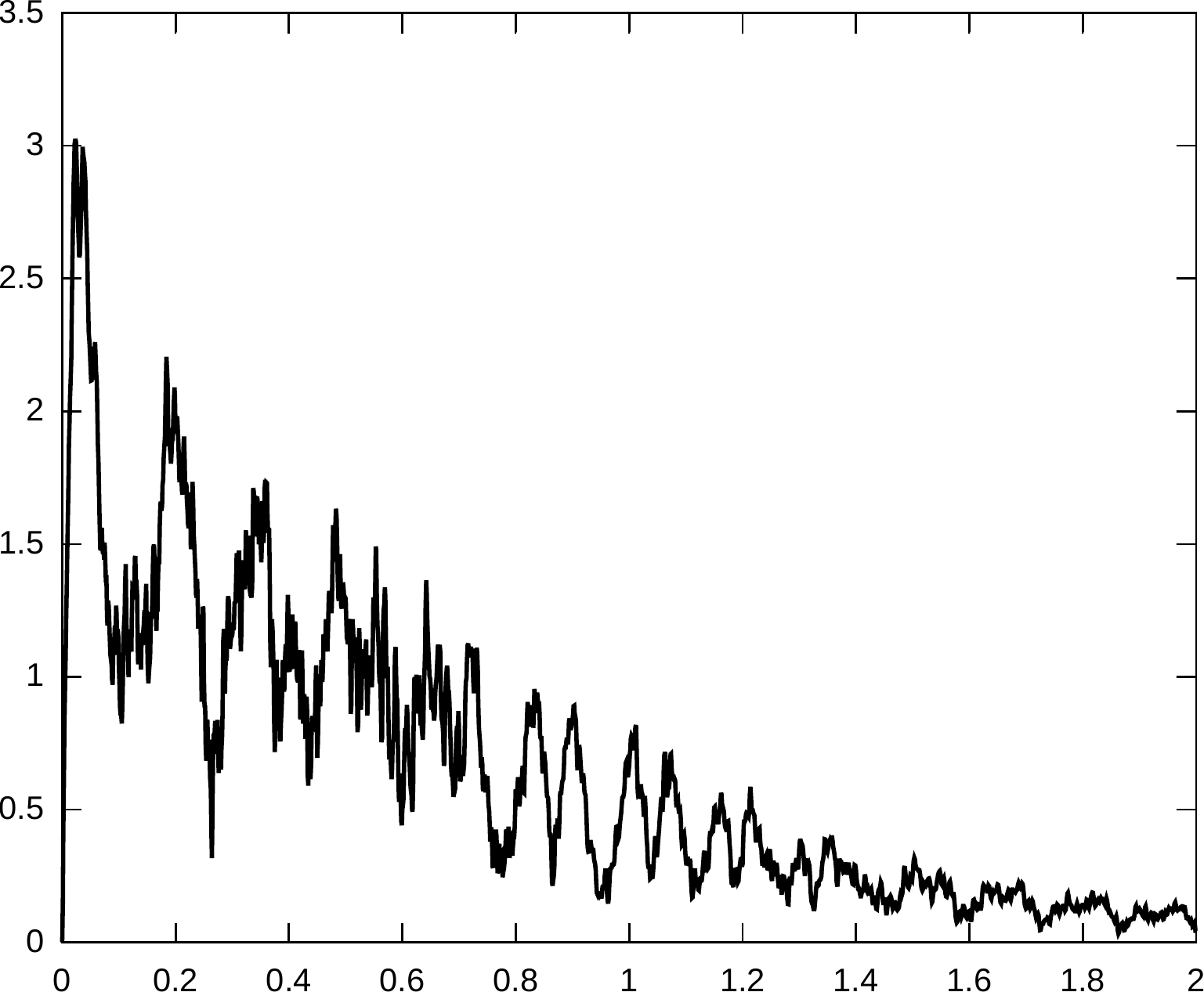}
   		\\\hline
   	\end{tabular}
   	\caption{Time evolution of target, control and kinetic energy  for the stationary (first row) and dynamic (second row) optimal control.}\label{energy}
   \end{figure}
Given a solution $q(s,t)$ of \eqref{sysintro} associated with either the static or dynamic optimal control, in Figure \ref{energy} we plot the time evolution of the three components of the cost functional $\mathcal J$: 
\begin{itemize}
    \item {\bf tip-target distance}: $\mathcal{J}_{q^*}(t):= |q(1,t)-q^\ast|^2. $
    
    We prioritized the reaching of the target by assigning an high value to the weight $\tau$ -- see Table \ref{parametertable}. As a consequence, in agreement with the theoretical setting, this is the component displaying the most prominent gap between the static and dynamic controls. We in particular see how, up to some oscillations, the dynamic control steers and keeps close to the target the tip $q(1,t)$ for most of the time of the evolution.
    
    \medskip
    \item  {\bf control energy}: $\mathcal{J}_u(t):=\int_0^1 u^2(s,t) ds\,.$
    
      In the case of static controls, $\mathcal{J}_{u^s}(t)$ is constant by construction. 
      The dynamic optimal controls are subject to the greatest variations in the beginning of the evolution and then stabilize around the static control energy, i.e., in agreement with the dissipative nature of the system, they actually converge to the optimal controls at the equilibrium. Finally note that dynamic controls perform better than static controls in terms of energy cost. In other words,    $\int_0^T\mathcal{J}_{u^d}(t)$ is smaller than   $\int_0^T\mathcal{J}_{u^s}(t)$.
      
     \medskip \item {\bf kinetic energy} $\mathcal{J}_v(t):=\frac{1}{2}\int_0^1\rho(s)|q_t(s,t)|^2 ds$. 
      Note that the evolutions of the $\mathcal J_v$ are comparable, but at final time $T$, dynamical optimal controls yield a kinetic energy $\mathcal{J}_v(T)$ slightly slower than the one associated to the stationary optimal controls. This is consistent with the fact that the cost functional $\mathcal J$ depends  on
   the kinetic energy only \emph{at final time}.

%   Finally note that, since the system is converging to an equilibrium due to frictional forces, $\mathcal{J}_v(t)\to 0$ as $t\to+\infty$: it is then reasonable to expect a larger time frame to improve the performances of the kinetic energy in both dynamic and stationary  cases.      
   \end{itemize}

\begin{remark}
    We remark that in our tests the mass distribution $\rho$ has an exponential decay, modeling a three-dimensional device with exponentially decaying thickness, see Section \ref{model}. This choice is motivated by the fact that  such structure can be viewed as an interpolation of self-similar, discrete hyper-redundant manipulators \cite{ems}, that is, robots composed by identical, rescaled modules. Besides the advantage in terms of design, self-similarity gives access to fractal geometry techniques \cite{lai}, allowing for a detailed investigation of inverse kinematics of the self-similarity manipulator, see for instance \cite{fibonacci}.
\end{remark}    

 \section{\uppercase{Optimal grasping}}\label{grasping}
 In this section we address a static optimal grasping problem, i.e., we look for stationary solutions of \eqref{sysintro} minimizing the distance from a target object and, at the same time, an integral quadratic cost on the associated curvature controls.
 
 More formally, let us denote by $\Omega_0$ an open subset of $\RR^2$ representing the object to be grasped, and by $\mbox{dist}(\cdot,\partial\Omega_0)$ the distance function from its boundary $\partial\Omega_0$. Moreover, we denote by $\chi_{\Omega_0}(\cdot)$ and $\chi_{\Omega_0^c}(\cdot)$, respectively the characteristic functions of $\Omega_0$ and its complement $\Omega_0^c$ in $\RR^2$. We consider the optimal control problem
\begin{equation}\label{contactfunctionalstationary}
\min \mathcal G,\quad\text{subject to \eqref{reducedstationary} and to $|u|\leq 1$.}
\end{equation} 
 where
 \begin{equation}\label{contactfunctional}
\begin{split}
    \mathcal G(q,u):=&\frac12\int_{[0,1]\setminus I} u^2 ds\\&+\frac{1}{2\tau}\int_0^1 
 \mbox{dist}^2(q(s),\partial\Omega_0)\left(\chi_{\Omega_0}(q(s))
 +\mu_0(s)\chi_{\Omega_0^c}(q(s))\right)ds\,
\end{split} ,\end{equation}
 is the cost functional for the grasping problem. The first integral term of $\mathcal G$ is a quadratic cost on the controls, computed only on the controlled region $[0,1]\setminus I$, see Section \ref{stat}. The second term of $\mathcal G$ aims to minimize the distance from the boundary of the target object \emph{without compenetration}. 
 %  , obtained by replacing the tip-target distance in the  functional \eqref{staticfunctional} by the following, more general term
%  \begin{equation}
% \label{contact}     
%  \frac{1}{2\tau}\int_0^1 
%  \mbox{dist}^2(q(s),\partial\Omega_0)\left(\chi_{\Omega_0}(q(s))+
%  \mu_0(s)\chi_{\Omega_0^c}(q(s))\right)ds\,.
%   \end{equation}
      More precisely, the first contribution, given by $\chi_{\Omega_0}$, acts as an obstacle, forcing all the points of the manipulator to move outside $\Omega_0$. The second term attracts points outside $\Omega_0$ on its boundary, according to $\mu_0$, a non negative weight describing which parts of the manipulator are preferred for grasping. In particular, if $\Omega_0$ is just a point and $\mu_0$ is a Dirac delta concentrated at $s=1$, we formally recover the tip-target distance.  On the other hand, different choices of $\mu_0$ allow one to obtain very different behaviors, as shown in the following tests -- see Table \ref{grasptable}.
%  From a numerical point of view, the optimization of \eqref{contact} requires the computation of the distance $\mbox{dist}(\cdot,\partial\Omega_0)$. It is well known that such a function admits an analytic expression only in some special cases, for instance if $\Omega_0$ is a circle or a square.

\subsubsection{Numerical tests for static optimal grasping.} 
\begin{figure}[!h]
 	\centering
 	\begin{tabular}{cc}
 		\includegraphics[width=.43\textwidth]{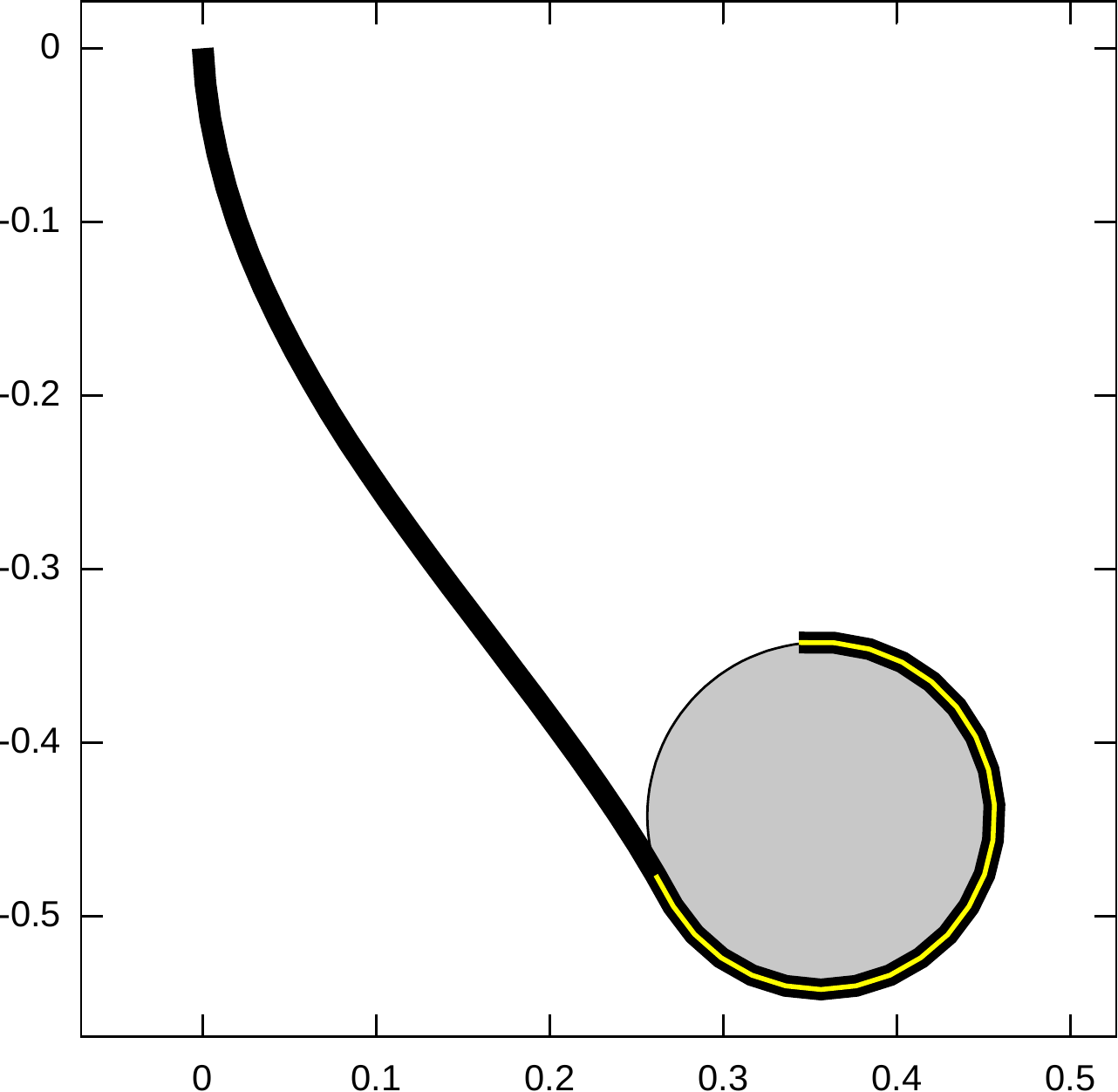}& 		\includegraphics[width=.523\textwidth]{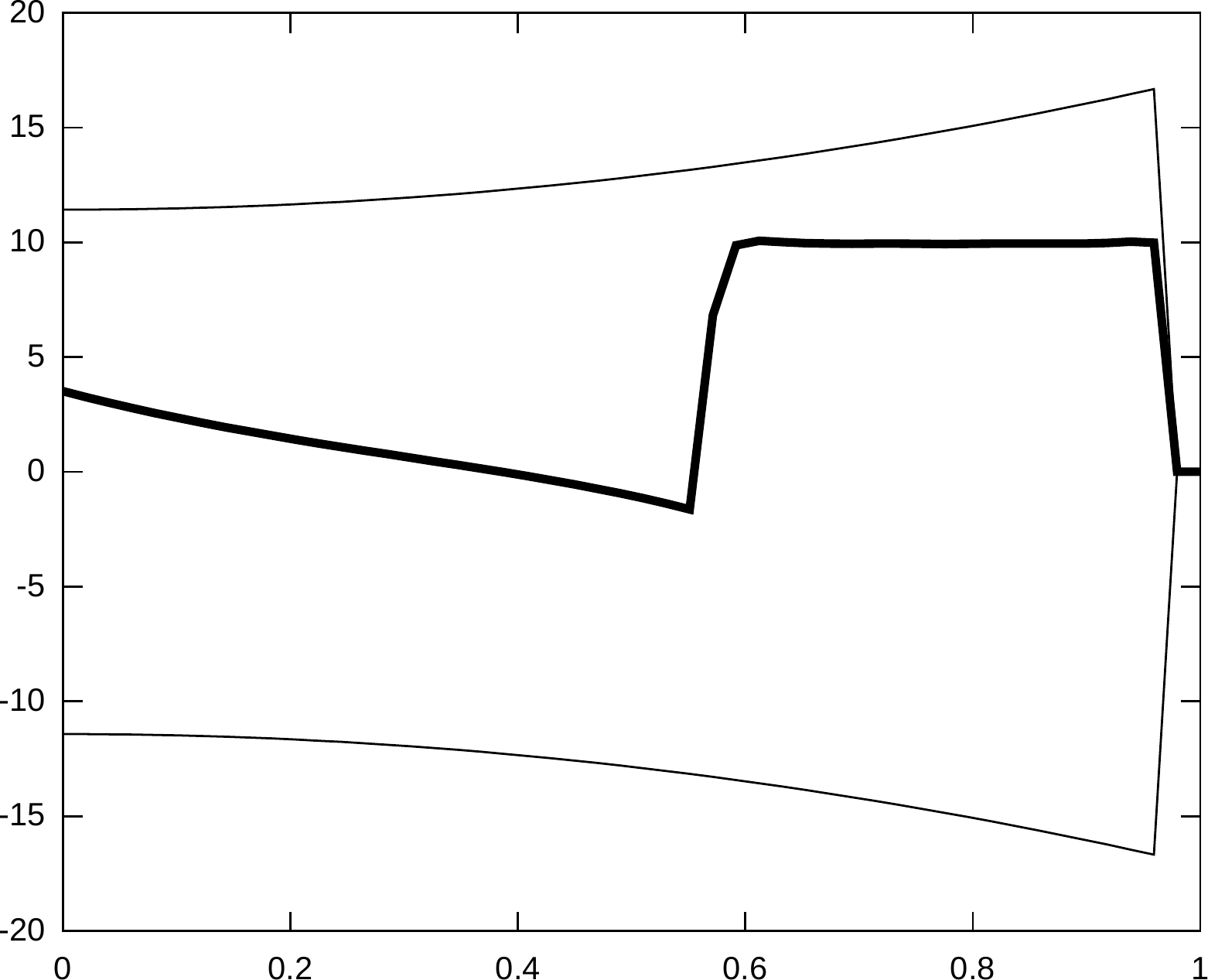}\\
 		(a)&	(b)	
 	\end{tabular}
 	\caption{\label{circle1} In (a) the solution $q$ of Test 5, in (b) the related signed curvature $\kappa(s)$.} 
 \end{figure}

We first consider the case in which $\Omega_0$ is a circle of radius $r_0$, and we set $\mu_0(s)=\chi_{[s_0,1]}(s)$ for some $0\le s_0 \le 1$. In this way, we expect the manipulator to surround the circle using only its terminal part of length $1-s_0$.
  \begin{table}[t]
  	\caption{Settings for the target $\Omega_0$ and for the grasping weight, $\mu_0$, where $s_0=0.55$ and $r_0=0.1$. The barycenter of $\Omega_0$ is $q^*=(0.3563,-0.4423)$, that is the target point for the tests in Section \ref{stat}. 
  	\label{grasptable}}
  \centering
  {	\renewcommand\arraystretch{1.3}\begin{tabular}{|c|c|c|}\hline
  		{\,\,\bf Test}& {\,\,\bf $\Omega_0$}&{\,\,\bf $\mu_0$}\\
  		\hline
  	\,\,Test 5\,\,& \,\,Circle of radius $r_0$\,\,&\,\,$\chi_{[s_0,1]}$\,\,\\
  	\,\,Test 6\,\,& \,\,Circle of radius $r_0$\,\,&\,\,$\delta_{s_0}+\delta_1$\,\,\\
  	\,\,Test 7\,\,& \,\,Square of side $2r_0$\,\,&\,\,$\chi_{[s_0,1]}$\,\,\\
  	\,\,Test 8\,\,& \,\,Square of side $2r_0$\,\,&\,\,$\delta_{s_0}+\delta_{\frac{s_0+1}{2}}+\delta_1$\,\,\\	\hline
  	\end{tabular}}
    \end{table}
Figure \ref{circle1} shows the results for Test 5, corresponding to the choice $r_0=0.1$ and $s_0=0.55$. Note that the length of the active part is $1-s_0=0.45$, namely less than the circumference $2\pi r_0\simeq 0.62$ of $\Omega_0$, hence the manipulator can not grasp the whole circle. Nevertheless, the curvature of  $\Omega_0$ is equal to $\frac{1}{r_0}=10$, which is below the upper bound $\bar\omega$ on the curvature $\kappa$ of the manipulator. This results in a good grasping, indeed we can observe a plateau in the graph of $\kappa$ right around the value $10$.
 
\begin{figure}[!h]
 	\centering
 	\begin{tabular}{cc}
 		\includegraphics[width=.43\textwidth]{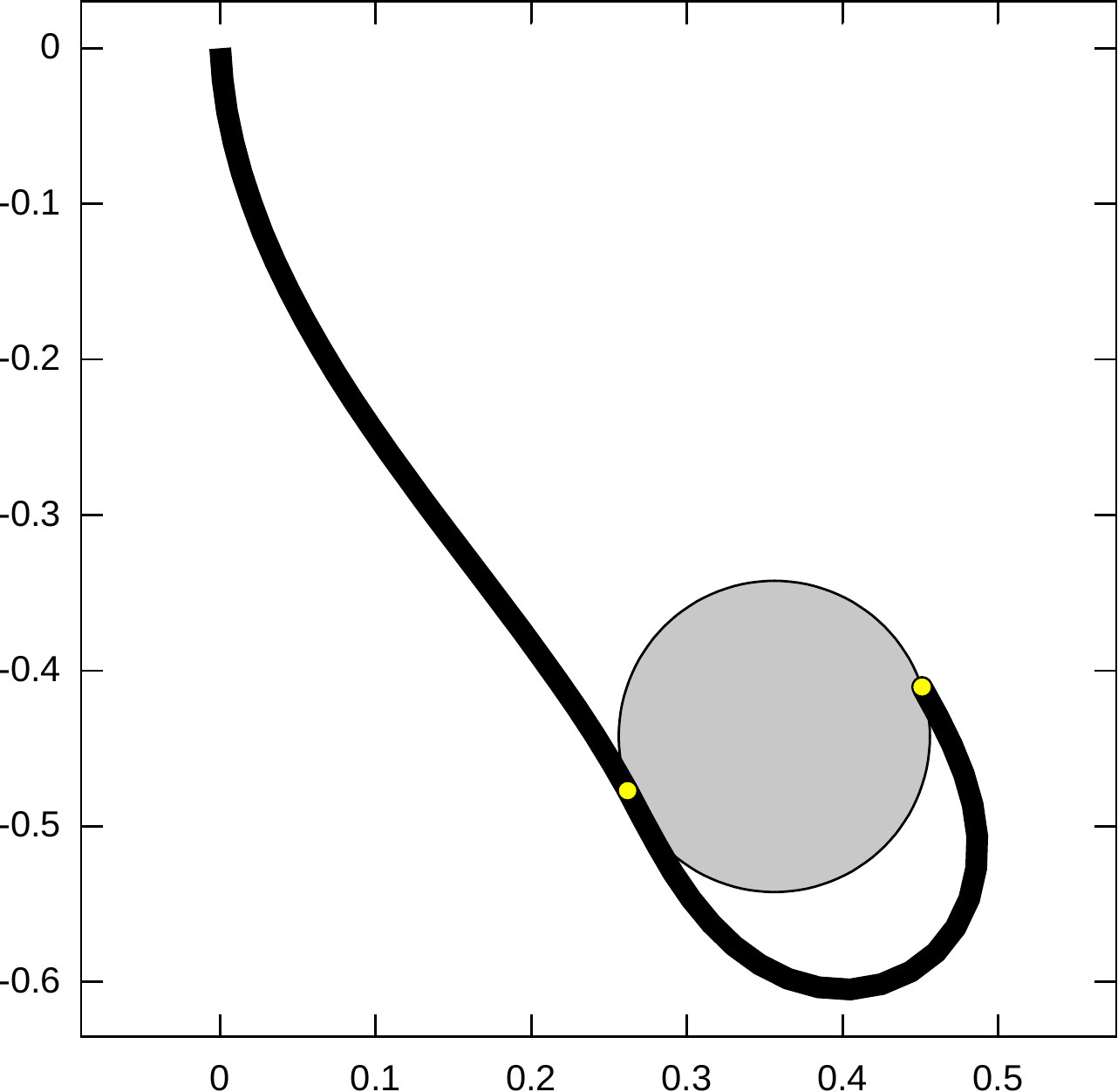}& 		\includegraphics[width=.523\textwidth]{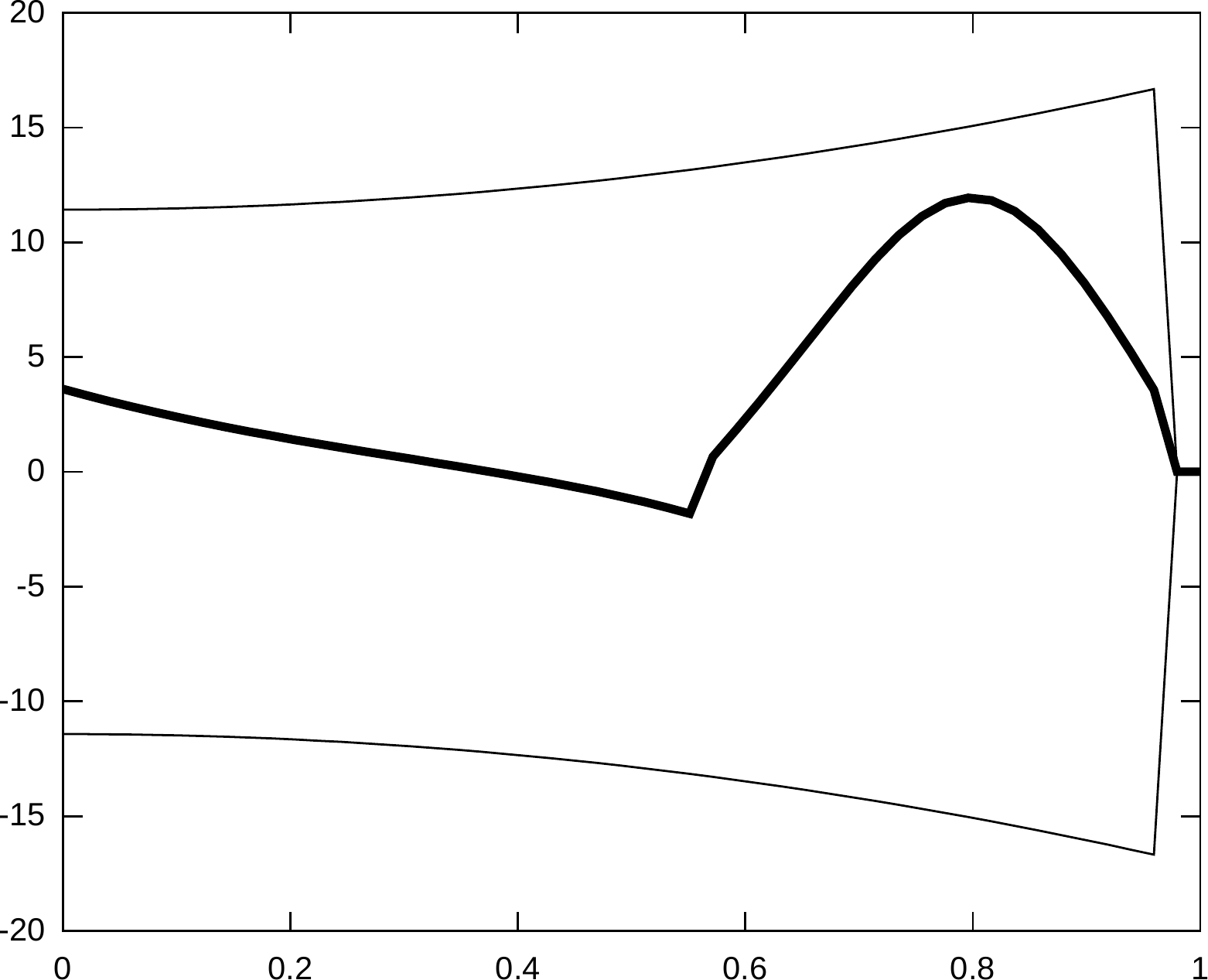}\\
 		(a)&	(b)	
 	\end{tabular}
 	\caption{\label{circle2} In (a) the solution $q$ of Test 6, in (b) the related signed curvature $\kappa(s)$.} 
 \end{figure}

In Test 6, we choose only two active points at the extrema of the interval $[s_0,1]$, namely we set $\mu_0(s)=\delta_{s_0}(s)+\delta_1(s)$. In Figure \ref{circle2}, we show the corresponding solution. Note that the maximum of $\kappa$ is now slightly larger than $10$, while the curvature energy is better optimized, since the arc between the active points no longer needs to be attached to $\Omega_0$. 

\begin{figure}[!h]
 	\centering
 	\begin{tabular}{cc}
 		\includegraphics[width=.43\textwidth]{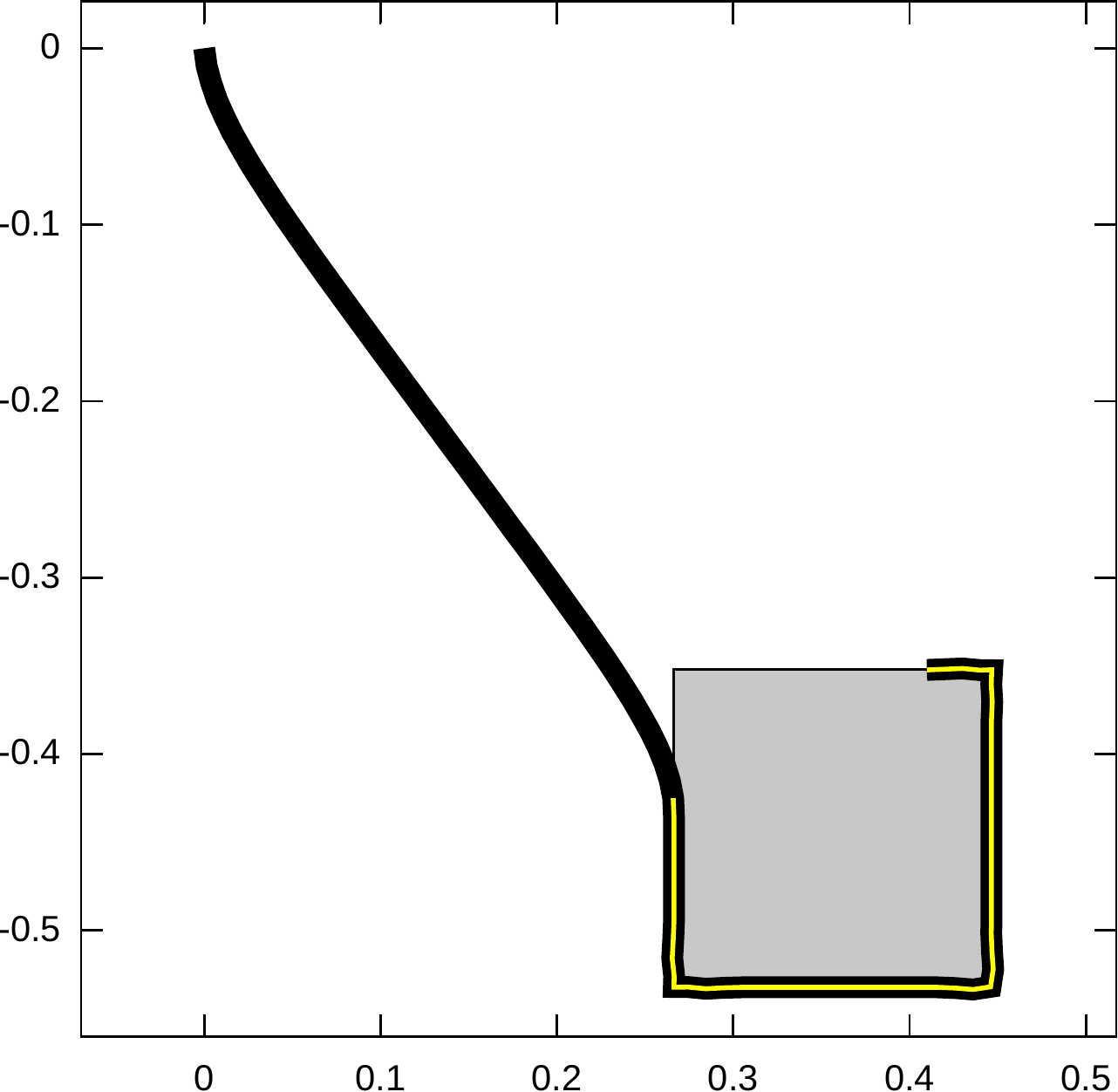}& 		\includegraphics[width=.523\textwidth]{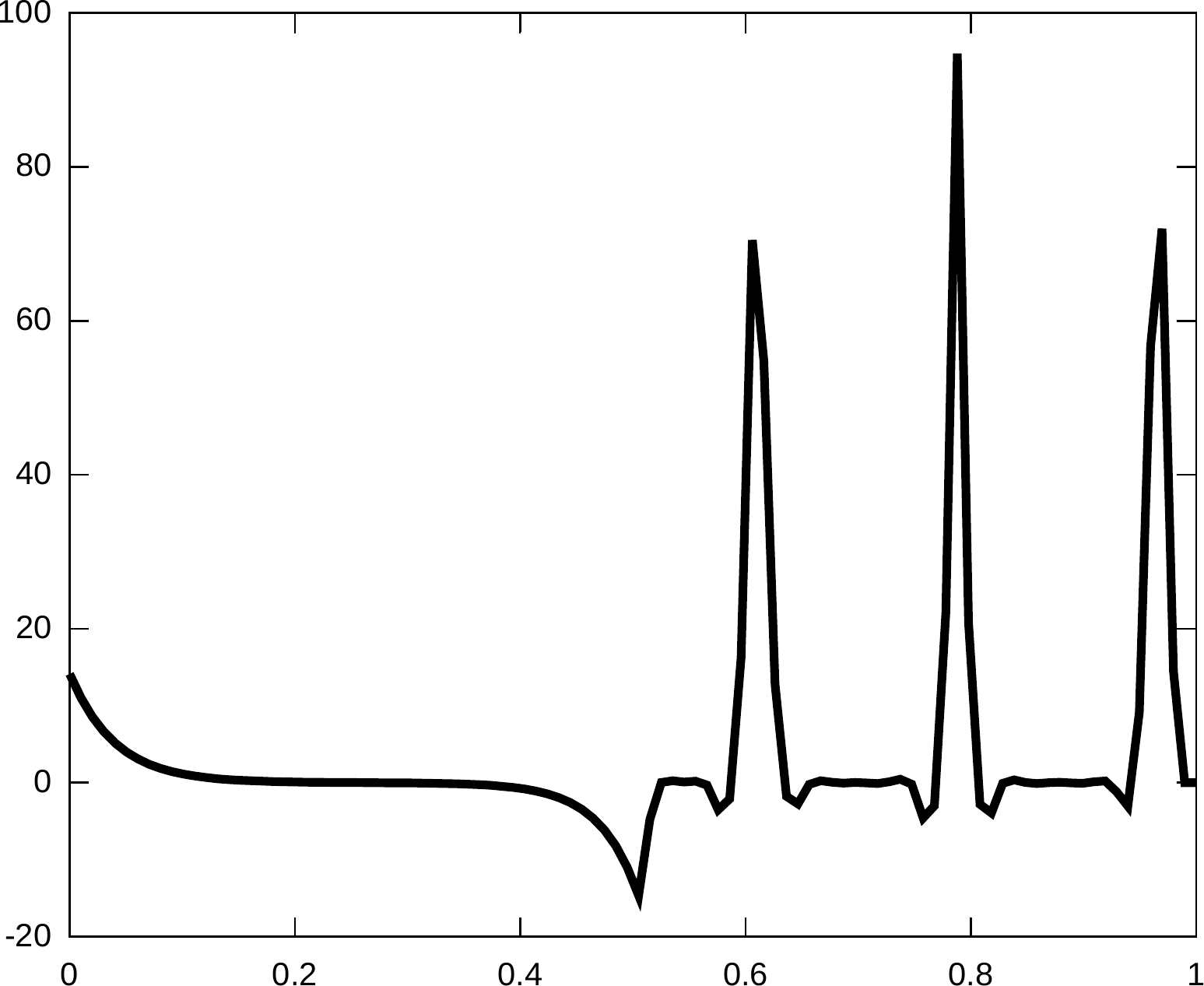}\\
 		(a)&	(b)	
 	\end{tabular}
 	\caption{\label{square1} In (a) the solution $q$ of Test 7, in (b) the related signed curvature $\kappa(s)$.} 
 \end{figure}

 We now consider Test 7, that is the case in which $\Omega_0$ is a square of side $2 r_0$, 
 and we set again $\mu_0(s)=\chi_{[s_0,1]}(s)$ for $s_0=0.55$ and $r_0=0.1$. This example is more challenging than the previous one, since a good grasping at the sharp corners of the square implies the divergence of the curvature of the manipulator. That is why we neglect the curvature constraints, reporting the results in Figure \ref{square1}. We clearly observe the presence of three spikes in the graph of $\kappa$, corresponding to the contact points at three corners. Moreover, we recognize some nearly flat parts of $\kappa$ corresponding to the sides of the square.  
 
\begin{figure}[!h]
 	\centering
 	\begin{tabular}{cc}
 		\includegraphics[width=.43\textwidth]{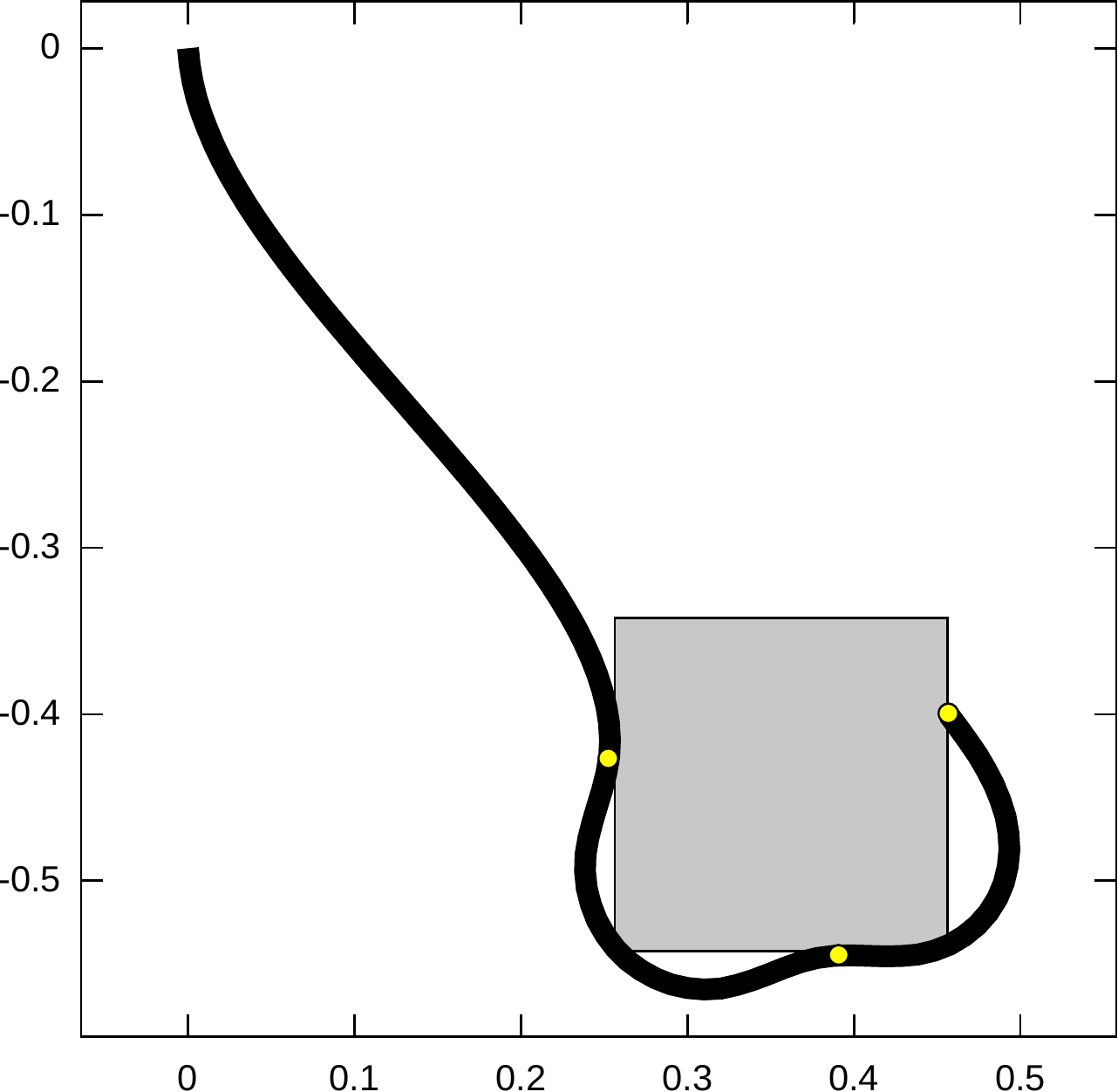}& 		\includegraphics[width=.523\textwidth]{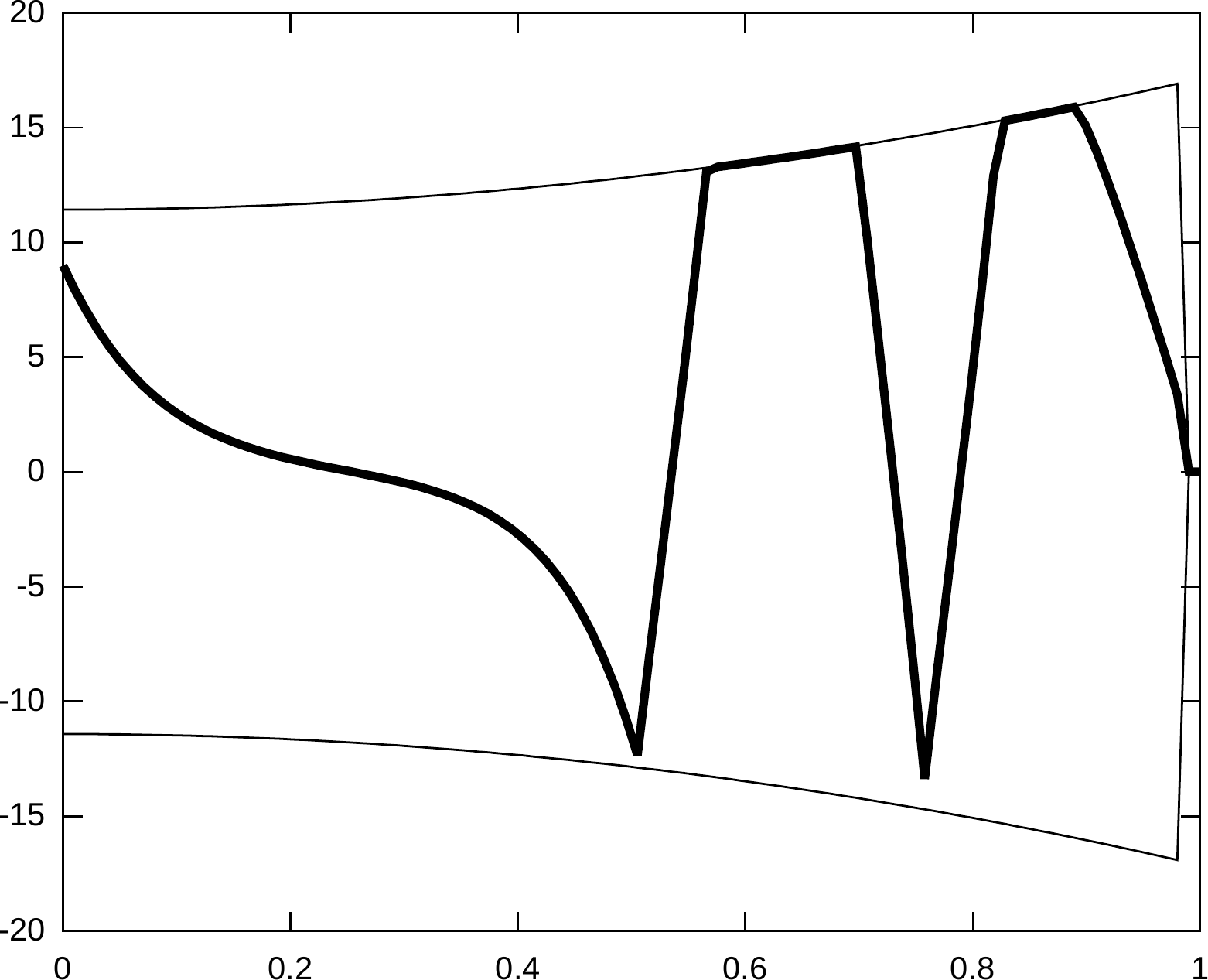}\\
 		(a)&	(b)	
 	\end{tabular}
 	\caption{\label{square2} In (a) the solution $q$ of Test 8, in (b) the related signed curvature $\kappa(s)$.} 
 \end{figure}

Finally, in Test 8 we choose only three active equispaced points in the interval $[s_0,1]$, namely setting $\mu_0(s)=\delta_{s_0}(s)+\
\delta_{\frac{s_0+1}{2}}(s)+\delta_1(s)$. Moreover, we restore the constraint on the  maximal curvature of the manipulator. In Figure \ref{square2}, we observe that the optimized configuration results from a non trivial balance between the contact energy and the curvature energy.

\section{\uppercase{Conclusions}}\label{conclusions}
In this paper we investigated a control model for the symmetry axis of a planar soft manipulator. The key features of the model (inextensibility, bending moment and curvature constraints and control) determine, together with internal and enviromental reaction forces, an evolution described by a system of fourth order nonlinear, evolutive PDEs. In particular, the equations of motion are derived as a formal limit of a discrete system, which in turn models a planar hyper-redundant manipulator subject to analogous, discrete, constraints. The comparison between the discrete (and rigid) and continuous (soft) model is a novelty in the investigation of the model, earlier introduced in \cite{CLL18}. We then addressed optimal reachability problems in both a static and dynamic setting: we tested the model in the case of uncontrolled regions and we showed that the model is also suitable for the investigation at the equilibrium of optimal inverse kinematics for hyper-redundant manipulators.
We then turned to optimal grasping problems, i.e., the problem of grasping an object (without compenetrating it) while minimizing a quadratic cost on the curvature controls. Our setting is a generalization of the earlier discussed reachability problem, and it allows for the investigation of optimal grasping for a quite large variety of objects. Moreover, the model allows to prioritize the preferred contact regions of the manipulator with the grasped object: we presented several numerical tests to illustrate this feature.

Future perspectives include stationary grasping problems for more complex target objects, possibly characterized by concavities and irregular boundaries. Our approach sets the ground for the search of optimal \emph{dynamic} controls: we plan to explore this in the future also with other optimization techniques, including model predictive control and machine learning algorithms.

 \bibliographystyle{splncs04}
 \bibliography{vurpo}

\begin{thebibliography}{10}
\providecommand{\url}[1]{\texttt{#1}}
\providecommand{\urlprefix}{URL }
\providecommand{\doi}[1]{https://doi.org/#1}

\bibitem{CLL18}
{Cacace}, S., {Lai}, A.C., {Loreti}, P.: Modeling and optimal control of an
  octopus tentacle. To appear on Siam Journal on Control and Optimization

\bibitem{CLL19}
Cacace, S., Lai, A.C., Loreti, P.: Control strategies for an octopus-like soft
  manipulator. In: Proceedings of the 16th International Conference on
  Informatics in Control, Automation and Robotics - Volume 1: ICINCO,. pp.
  82--90. INSTICC, SciTePress (2019). \doi{10.5220/0007921700820090}

\bibitem{hyper}
Chirikjian, G.S., Burdick, J.W.: An obstacle avoidance algorithm for
  hyper-redundant manipulators. In: Proceedings., IEEE International Conference
  on Robotics and Automation. pp. 625--631. IEEE (1990)

\bibitem{hyper2}
Chirikjian, G.S., Burdick, J.W.: The kinematics of hyper-redundant robot
  locomotion. IEEE transactions on robotics and automation  \textbf{11}(6),
  781--793 (1995)

\bibitem{kinematicssoft}
Jones, B.A., Walker, I.D.: Kinematics for multisection continuum robots. IEEE
  Transactions on Robotics  \textbf{22}(1),  43--55 (2006)

\bibitem{dynamicoctopusrobot}
Kang, R., Kazakidi, A., Guglielmino, E., Branson, D.T., Tsakiris, D.P.,
  Ekaterinaris, J.A., Caldwell, D.G.: Dynamic model of a hyper-redundant,
  octopus-like manipulator for underwater applications. In: Intelligent Robots
  and Systems (IROS), 2011 IEEE/RSJ International Conference on. pp.
  4054--4059. IEEE (2011)

\bibitem{octopusrobotprescribed}
Kazakidi, A., Tsakiris, D.P., Angelidis, D., Sotiropoulos, F., Ekaterinaris,
  J.A.: Cfd study of aquatic thrust generation by an octopus-like arm under
  intense prescribed deformations  \textbf{115},  54--65 (2015)

\bibitem{lai}
Lai, A.C.: Geometrical aspects of expansions in complex bases. Acta Mathematica
  Hungarica  \textbf{136}(4),  275--300 (2012)

\bibitem{hand}
Lai, A.C., Loreti, P.: {Robot's hand and expansions in non-integer bases}.
  {Discrete Mathematics \& Theoretical Computer Science}  \textbf{16}(1) (Jun
  2014), \url{https://dmtcs.episciences.org/3913}

\bibitem{icinco1}
Lai, A.C., Loreti, P.: Discrete asymptotic reachability via expansions in
  non-integer bases. In: 2012 9-th international conference on informatics in
  control, automation and robotics (ICINCO). vol.~2, pp. 360--365. IEEE (2012)

\bibitem{icinco2}
Lai, A.C., Loreti, P., Vellucci, P.: A model for robotic hand based on
  fibonacci sequence. In: 2014 11th International Conference on Informatics in
  Control, Automation and Robotics (ICINCO). vol.~2, pp. 577--584. IEEE (2014)

\bibitem{ems}
Lai, A.C., Loreti, P., Vellucci, P.: A continuous fibonacci model for robotic
  octopus arm. In: 2016 European Modelling Symposium (EMS). pp. 99--103. IEEE
  (2016)

\bibitem{fibonacci}
Lai, A.C., Loreti, P., Vellucci, P.: A fibonacci control system with
  application to hyper-redundant manipulators. Mathematics of Control, Signals,
  and Systems  \textbf{28}(2), ~15 (2016)

\bibitem{optimalpdebook}
Tr{\"o}ltzsch, F.: Optimal control of partial differential equations: theory,
  methods, and applications, vol.~112. American Mathematical Soc. (2010)

\end{thebibliography}

\end{document}